\documentclass[12pt,a4paper,draft]{article}

\usepackage[english]{babel}
\usepackage[psamsfonts]{amssymb}
\usepackage{cmmib57}
\usepackage{exscale}
\usepackage{amsmath}
\usepackage{amscd}
\usepackage{latexsym}
\usepackage{amsthm}
\usepackage{amssymb}
\usepackage[dvips]{graphicx}

\numberwithin{equation}{section}

\newcommand{\R}{\mathbb R}
\newcommand{\N}{\mathbb N}

\newcommand{\cb}{\mathcal B}
\newcommand{\cac}{\mathcal{C} }

\newcommand{\ce}{\mathcal E}
\newcommand{\cf}{\mathcal F}

\newcommand{\ci}{\mathcal I}
\newcommand{\cj}{\mathcal J}

\newcommand{\cs}{\mathcal S}
\newcommand{\cw}{\mathcal W}

\newcommand{\ck}{\mathcal K}
\newcommand{\cp}{\mathcal P}

\newcommand{\vp}{\varphi}
\newcommand{\vt}{\vartheta}

\newcommand{\al}{\alpha}
\newcommand{\Om}{\Omega}
\newcommand{\om}{\omega}

\newcommand{\Del}{\Delta}

\newcommand{\lam}{\lambda}
\newcommand{\Lam}{\Lambda}
\newcommand{\Gam}{\Gamma}

\newcommand{\ffi}{\varphi}

\newcommand{\cinf}{\mathcal{C}^\infty}
\newcommand{\IN}{\mathbb{N}}
\newcommand{\re}{\mathbb{R}}
\newcommand{\red}{\mathbb{R}^{d}}
\newcommand{\sig}{\sigma}

\newcommand{\punt}{\cdot}
\newcommand{\tf}{\mathcal{F}}
\newcommand{\beq}{\begin{equation}}
\newcommand{\eeq}{\end{equation}}

\newtheorem{theorem}{Theorem}[section]

\newtheorem{definition}[theorem]{Definition}

\newtheorem{hypothesis}[theorem]{Hypothesis}
\newtheorem{lemma}[theorem]{Lemma}

\newtheorem{proposition}[theorem]{Proposition}

\newtheorem{remark}[theorem]{Remark}

\title{\bf Stochastic integrals for spde's: a comparison}
\author{{\bf Robert C.~Dalang}\thanks{Partly supported by the Swiss National Foundation for Scientific Research.}\\ Ecole Polytechnique F\'ed\'erale de Lausanne \and {\bf Llu\'is Quer-Sardanyons}\thanks{Partly supported by the grant MEC-FEDER Ref. MTM20006-06427 from the
  Direcci\'on General de Investigaci\'on, Ministerio de Educaci\'on y Ciencia, Spain.}\\ Universitat Aut\`onoma de Barcelona}
\date{}
\begin{document}
\maketitle
\abstract{We present the Walsh theory of stochastic integrals with respect to martingale measures, alongside of the Da Prato and Zabczyk theory of stochastic integrals with respect to Hilbert-space-valued Wiener processes and some other approaches to stochastic integration, and we explore the links between these theories. We then show how each theory can be used to study stochastic partial differential equations, with an emphasis on the stochastic heat and wave equations driven by spatially homogeneous Gaussian noise that is white in time. We compare the solutions produced by the different theories.}
\vskip 1in

\eject

\tableofcontents

\eject


\section{Introduction}

   The theory of stochastic partial differential equations (spde's) developed on the one hand, from the work of J.B.~Walsh \cite{walsh}, and on the other hand, through work on stochastic evolution equations in Hilbert spaces, such as in \cite{dawson}. An important milestone in the latter approach is the book of Da Prato and Zabczyk \cite{dz}. 
   
    These two approaches led to the development of two distinct schools of study for spde's, based on different theories of stochastic integration: the Walsh theory, which emphasizes integration with respect to worthy martingale measures, and a theory of integration with respect to Hilbert-space-valued processes, as expounded in \cite{dz}. A consequence of the presence of these separate theories is that the literature published by each of the two schools is difficult to access when one has been trained in the other school. This is unfortunate since both approaches have advantages and in some problems, using both approaches can be useful (one example of this is \cite{dmz}).
    
    The objective of this paper is to help create links between these two schools of study. It is addressed to researchers who have some familiarity with at least one of the two approaches. We develop both theories, and explore the links between the two. Then we show how each theory is used to study spde's. The Walsh theory emphasizes solutions that are random fields, while \cite{dz} centers around solutions in Hilbert spaces of functions. Each theory is presented rather succinctly, the main focus being on relationships between the theories. We show that these theories often (but not always) lead to the same solutions to various spde's.

   It should be mentioned that the general theory of integration with respect to Hilbert-space-valued processes and its generalizations---such as the stochastic integral with respect to cylindrical processes---was well-developed several years before \cite{walsh} and more than
a decade before reference \cite{dz} appeared: see, for instance, the book of M\'etivier and Pellaumail \cite{metivier}. This reference, and several others, are cited in \cite{dz} and \cite{walsh}. However, J.B.~Walsh preferred to develop his own integral, even though he realized that the two were related (see the Notes at the end of \cite{walsh}).

   Here, we present in Section \ref{integral} a modern formulation of the theory of stochastic integrals with respect to cylindrical Wiener processes, as developed in \cite{metivier}, as a unifying integral behind most of those that were introduced later on.  This integral is briefly recalled in Section \ref{integral}. In Section \ref{prel}, we show how spatially homogeneous noise that is white in time can be viewed as a cylindrical Wiener process on a particular Hilbert space. Emphasizing this type of noise is natural, since in recent years, following in particular the papers of Mueller \cite{mueller}, Dalang and Frangos \cite{DF}, Dalang \cite{dalang} and Peszat and Zabczyk \cite{pz1,pz2}, this type of noise has been used by several researchers. This is due in part to the fact that it leads to a theory of non-linear spde's in spatial dimensions greater than $1$, while non-linear spde's driven by space-time white noise generally only have a solution in spatial dimension $1$. In Section \ref{homo}, we show (Proposition \ref{w-d}) that the Walsh stochastic integral and the extension presented by Dalang \cite{dalang} and Nualart and Quer-Sardanyons \cite{nualartquer} can be viewed as integrals as defined in Section \ref{integral}. Section \ref{distri} gives a wide class of integrable processes. In Section \ref{ext-dm}, we  discuss the relationship between this integral and the function-valued stochastic integral introduced by Dalang and Mueller in \cite{dm}.
   

   In Section \ref{infdim}, we sketch the construction of the infinite dimensional stochastic integral in the setup of Da Prato and Zabczyk \cite{dz}. We also make use of the more recent presentation of Pr\'ev\^ot and R\"ockner \cite{pr}. In Section \ref{recall}, we recall some basic properties of Hilbert-Schmidt operators. Section \ref{wiener} gives the relationship between a Hilbert-space-valued Wiener process and a cylindrical Brownian motion, in the case where the covariance operator has finite trace. Hilbert-space-valued stochastic integrals are defined in Section \ref{daprato}. In particular, we show in Proposition \ref{prop3} how this infinite-dimensional stochastic integral can be written as a series of It\^o stochastic integrals. This is used in Section \ref{cylinf} to show how the integrals of Section \ref{big} can be interpreted in the infinite-dimensional context. The case of covariance operators with infinite-trace is discussed in Section \ref{cylwin}.
   
   It is well-known that in certain cases, the Hilbert-space-valued integral is equivalent to a martingale-measure stochastic integral. For instance, it is pointed out in \cite[Section 4.3]{dz} that when the random perturbation is space-time white noise, then Walsh's stochastic integral in \cite{walsh} is equivalent to an infinite-dimensional stochastic integral as in \cite{dz} (see also \cite{K}). Of course, space-time white noise is only a special case of spatially homogeneous noise, and we are interested in comparing solutions to spde's driven by this more general noise. The function-valued approach of \cite{dm} gives solutions to spde's for which it is not known if a random field solution exists, and the Hilbert-space approach is even more general. However, for a wide class of spde's that have solutions in two or more of these formulations, such as the stochastic heat equation ($d \ge 1$) and wave equation ($d \in \{1,2,3\}$) driven by spatially homogeneous noise, we will show that the solutions turn out to be equivalent. One does not expect this to be the case in all situations. Indeed, there are a few cases in which a solution exists with one approach and is known {\em not to exist} in one of the others. For instance, for noise concentrated on a hyperplane, as considered in \cite{DL},
the authors establish existence of function-valued solutions and show that there is no random field solution.

   In Section \ref{spde}, we first discuss the random field approach to the study of spde's, with an emphasis on the stochastic heat and wave equations. We use the stochastic integral of Section \ref{homo} to extend the result of \cite{dalang} to arbitrary initial conditions (Theorem \ref{existence}). We then discuss the Hilbert-space-valued approach to the study of the same equations, taking the stochastic wave equation and the approach of \cite{pz2} as a primary example. In particular, we show that the mild random field solution of Theorem \ref{existence}, when interpreted as a Hilbert-space-valued process, yields the solution given in \cite{pz2}. This is achieved by identifying the multiplicative non-linearity with an appropriate Hilbert-Schmidt operator, and using the relationships between stochastic integrals exposed in Section \ref{infdim}. Since the two solutions are defined using different Hilbert spaces, the embedding from one Hilbert space to the other has to be written explicitly. Finally, in Section \ref{relation-dm}, we compare the random field solution of the stochastic wave equation with the function-valued solution constructed in \cite{dm}. Again, in cases where both types of solutions are defined, that is, in spatial dimensions $d \in \{1,2,3\}$, we show that the random field solution yields the function-valued solution (Theorem \ref{real-dm}). Overall, Section \ref{spde} unifies the existing literature on the stochastic heat and wave equations driven by spatially homogeneous noise, and clarifies the relationships between the various approaches.

\section{Stochastic integrals with respect to a Gaussian spatially homogeneous noise}
\label{big}

In this section, we recall the notion of cylindrical Wiener process and the stochastic  
integral with  respect to such processes. Then we introduce a spatially homogeneous
Gaussian noise that is white in time, and show how to interpret this noise as a
cylindrical Wiener process.  Building on material presented in  \cite{nualartquer}, we
then relate the stochastic integral with respect to this particular cylindrical Wiener
process with Walsh's martingale measure stochastic integral and the extension given by
Dalang in \cite{dalang}. We also discuss the extension given in Dalang and Mueller \cite{dm}.

\subsection{Stochastic integration with respect to a cylindrical
Wie\-ner process} \label{integral}

Fix a Hilbert space $V$ with inner product $\langle \cdot,\cdot\rangle_V$. Following \cite{pardoux,metivier}, we define the general notion of cylindrical Wiener process in $V$.

\begin{definition}\label{cyl}
Let $Q$ be a symmetric (self-adjoint) and non-negative definite bounded linear
operator on $V$. A family of random variables $B=\{B_t(h),\, t\geq 0,\, h\in V\}$ is a \emph{cylindrical Wiener process} on $V$ if the following
two conditions are fulfilled:

  1.~for any $h\in V$, $\{B_t(h),\ t\geq 0\}$ defines a Brownian motion
with variance $t\langle Qh,h\rangle_V$;

   2.~for all $s,t\in \re_+$ and $h,g\in V$,
$$
   E\left( B_s(h) B_t(g)\right)=(s\wedge t)\langle Qh,g\rangle_V,
$$
where $s\wedge t := \min(s,t)$. If $Q=\mbox{\rm Id}_V$ is the identity operator in $V$, then $B$ will be called a
\emph{standard cylindrical Wiener process}. We will refer to $Q$ as the {\em covariance} of $B$.
\end{definition}

Let $\mathcal{F}_t$ be the $\sigma$-field generated by the random
variables $\{B_s(h),\, h\in V,\, 0\le s\le t\}$ and the
$P$-null sets. We define the predictable $\sigma$-field as the
$\sigma$-field in $[0,T]\times \Om$ generated by
the sets $\{ (s,t]\times A,\; A\in \mathcal{F}_s,\; 0\le s<t\le T \}$.

We denote by $V_Q$ the Hilbert space $V$ endowed with the
inner-product
$$\langle h,g\rangle_{V_Q}:=\langle Q h,g\rangle_V,\qquad h,g\in V.$$
We can now define the stochastic integral of any predictable square-integrable process
with values in $V_Q$, as follows. Let $(v_j)_j$ be a complete orthonormal basis of the
Hilbert space $V_Q$. For any predictable process $g\in L^{2}(\Omega \times [0,T];V_Q)$,
it turns out that the following series is convergent in $L^2(\Om, \tf,P)$ and the sum
does not depend on the chosen orthonormal system:
\beq
   g\cdot B :=
\sum_{j=1}^\infty \int_0^T \langle g_s ,v_j\rangle_{V_Q}\, dB_s(v_j). \label{22}
\eeq
We notice that each summand in the above series is a classical It\^o integral with respect
to a standard Brownian motion, and the resulting stochastic integral is a real-valued
random variable. The stochastic integral $g \cdot B$
is also denoted by $\int_0^T g_s\, dB_s$. The independence of the terms in the series \eqref{22} leads to the isometry property
$$
E\left( \left( g\cdot B\right)^{2}\right) = E\left( \left(\int_0^T
g_s\, dB_s\right)^2 \right)=  E\left( \int_{0}^{T}\left\|
g_s\right\| _{V_Q}^{2}\, ds\right) .
$$

We note that there is an alternative way of defining this integral: one can start by defining the
stochastic integral in (\ref{22}) for a class of \emph{simple}
predictable $V_Q$-valued processes, and then use the isometry property to extend the integral to elements of $L^2(\Om\times[0,T]; V_Q)$ by checking that these simple processes are dense in this set.

\subsection{Spatially homogeneous noise as a cylindrical Wiener process}
\label{prel}

We now define the Gaussian random noise that will play a central role in this paper. On a complete probability space $(\Omega, \mathcal{F},P)$, we consider a family of mean zero Gaussian random variables
$W = \{W(\varphi),\, \varphi \in \mathcal{C}_0^\infty (\mathbb{R}^{d+1})\}$, where $\mathcal{C}_0^\infty (\mathbb{R}^{d+1})$ denotes the space of infinitely differentiable functions with compact support,
with covariance
\begin{align}\label{1}
E(W(\vp) W(\psi)) &= \int_0^\infty dt \int_{\red} dx \int_{\red} dy\, \vp(t,x) f(x-y) \psi(t,y)\\
 &= \int_0^\infty dt \int_{\red} dx\,  f(x)\, (\vp(t)\ast \tilde\psi(t))(x),
\nonumber
\end{align}
where ``$\ast$" denotes convolution in the spatial variable and $\tilde \vp (x):= \vp(-x)$.

In the above, $f$ is a non-negative and non-negative definite continuous function on $\mathbb{R}^{d}\setminus \{0\}$ which is integrable in a neighborhood of $0$ and is the Fourier
transform of a non-negative tempered measure $\mu$ on $\mathbb{R}^{d}$. That is, by definition of the Fourier transform on the space $\cs'(\red)$ of tempered distributions, for all $\ffi$ belonging to the space $\cs(\red)$ of rapidly decreasing $\cinf$ functions,
$$\int_{\red} f(x) \ffi(x)\, dx=\int_{\red} \tf \ffi(\xi)\, \mu(d\xi),$$
and there is an integer $m\geq 1$ such that
\begin{equation}\label{tempered}
\int_{\red} (1+|\xi|^2)^{-m}\mu(d\xi) <\infty.
\end{equation}
We have denoted by $\tf \ffi$ the Fourier transform of $\ffi\in
\cs(\red)$:
$$\tf \ffi(\xi)= \int_{\red} \ffi(x) e^{-2\pi i\xi\cdot x}
\, dx.$$ The measure $\mu$ is called the {\em spectral measure} of $W$ and is necessarily symmetric (see \cite[Chap. VII, Th\'eor\`eme XVII]{schwartz}), and $f$ is necessarily an even function.
The covariance (\ref{1}) can also be written, using elementary
properties of the Fourier transform, as
$$
  E(W(\varphi) W(\psi)) = \int_0^\infty dt \int_{\mathbb{R}^{d}}  \mu(d\xi)\, \mathcal{F} \varphi(t)(\xi)\, \overline{\mathcal{F} \psi(t)(\xi)} .
$$

\begin{remark}
 We observe that, in (\ref{1}), it is possible to take a slightly more general spatial correlation: the function $f$ could be replaced by a non-negative and non-negative definite tempered measure (for instance, see \cite[Section 2]{dm}). Formula \eqref{1}, which we use for the sake of clarity in the exposition, corresponds to the case where this measure is absolutely continuous with respect to Lebesgue measure on $\red$, with density $f$. 
\label{rem0}
\end{remark}

  It is natural to associate a Hilbert space with $W$: let $U$ the completion of the Schwartz space $\cs (\red)$
endowed with the semi-inner product
\begin{equation}\label{ip}
\langle \vp,\psi \rangle_{U} =\int_{\red} dx  \int_{\red} dy\, \vp(x) f(x-y) \psi(y) =\int_{\red} \mu(d\xi)\, \tf \vp(\xi) \overline{\tf \psi(\xi)},
\end{equation}
$\vp,\psi\in \cs (\red)$, and associated semi-norm $\Vert \cdot
\Vert_U $.
Then $U$ is a Hilbert space that may contain Schwartz distributions (see \cite[Example 6]{dalang}).

\begin{remark}
Let $\tilde L^2(\red,d\mu)$ be the subspace of $L^2(\red,d\mu)$ consisting of functions $\phi$ such that $\tilde \phi=\phi$. It is not difficult to check that one can identify $U$ with the set $\{\Psi \in \cs'(\red): \Psi = \tf^{-1}\phi,\mbox{ where } \phi\in \tilde L^2(\red,d\mu)\}$, with inner
product
$$\langle \tf^{-1} \phi,\tf^{-1} \ffi\rangle_U= \langle \phi,\ffi\rangle_{L^2(\red, d\mu)},\qquad \phi,\ffi\in \tilde L^2(\red, d\mu).$$
\label{rem1}
\end{remark}

We fix a time interval $[0,T]$ and we set $U_T:=L^2([0,T];U)$. This set is equipped with the norm given by
$$
   \Vert g \Vert^2_{U_T} = \int_0^T \Vert g(s) \Vert^2_U\,ds .
$$


   We now associate a cylindrical Wiener process to $W$, as follows. A direct calculation using \eqref{1} shows that the generalized
Gaussian random field $\{W(\vp),\, \vp\in \cinf_0([0,T]\times
\red)\}$ is a random linear functional, in the sense that $W(a \vp+b
\psi) = a W(\vp) + bW(\psi)$, and $\vp \mapsto W(\vp)$ is an
isometry from $(\cinf_0 ([0,T]\times \red), \Vert \cdot
\Vert_{U_T})$ into $L^2(\Omega,\tf,P)$. The following lemma
identifies the completion of $\cinf_0 ([0,T]\times \red)$ with
respect to $\Vert \cdot \Vert_{U_T}$.

\begin{lemma}\label{lem2.2}
The space $\cinf_0 ([0,T]\times \red)$ is dense in $U_T = L^2([0,T];U)$
for $\Vert \cdot \Vert_{U_T}$.
\end{lemma}

\noindent \emph{Proof.} Following \cite{nualartquer}, we will use the notation $\varphi_1(\cdot)$ to indicate that $\varphi_1$ is a function $t \mapsto \varphi_1(t)$ of the time-variable, and $\varphi_2(\star)$ to indicate that $\varphi_2$ is a function $x \mapsto \varphi_2(x)$ of the spatial variable.

   Let $C$ denote the closure of $\cinf_0 ([0,T]\times
\red)$ in $U_T$ for $\Vert \cdot \Vert_{U_T}$. Clearly, $C$ is a
subspace of $U_T$. The proof can be split into three parts.\\

\noindent \emph{Step 1.} We show that elements of $U_T$ of the form
$\vp_1(\cdot) \vp_2(\star )$, where $\vp_1\in \cinf_0(\re_+; \re)$ with support
included in $[0,T]$ and $\vp_2\in \cs (\red)$, belong to $C$. Using the fact that
$$
 \int_\red dx\, f(x) \, (\vert \vp_2 \vert  \ast \vert  \tilde \vp_2\vert )(x) < \infty
$$
because $f$ is a tempered function by \eqref{tempered} and $\vert \vp_2 \vert  \ast \vert  \tilde \vp_2\vert$ decreases rapidly, together with dominated convergence, one checks
that there is a sequence $(\vp_2^n)_n \subset \cinf_0(\red)$ such that
$\lim_{n\rightarrow \infty} \|\vp_2 - \vp_2^n\|_{U}=0$. Then, by the
very definition of the norm in $U_T$, one easily proves that
$\lim_{n\rightarrow \infty} \|\vp_1 \vp_2 - \vp_1
\vp_2^n\|_{U_T}=0$. Therefore, $\vp_1(\cdot) \vp_2(\star )\in U_T$.
\\

\noindent \emph{Step 2.} Suppose that we are given
$\vp_1\in L^2([0,T]; \re)$ and $\vp_2\in \cs (\red)$. We show that $\vp_1(\cdot)\vp_2(\star )\in C$. Indeed, let $(\vp_1^n)_n\in \cinf_0(\re_+)$ be such that,
for all $n$, the support of $\vp_1^n$ is contained in $[0,T]$ and
$\vp_1^n \rightarrow \vp_1$ in $L^2([0,T];\re)$. Then $\vp_1^n
\vp_2\in C$ by Step 1, and one checks that $\vp_1^n \vp_2$
converges, as $n$ tends to infinity, to $\vp_1 \vp_2$ in $U_T$.
Therefore, $\vp_1(\cdot)\vp_2(\star )\in C$.\\

\noindent \emph{Step 3.} Suppose that $\vp\in
U_T$. We show that $\vp\in C$. Indeed, let $(e_j)_j$
be a complete orthonormal basis of $U$ with $e_j\in \cs(\red)$, for
all $j$. Then, since $\vp(s)\in U$ for any $s\in [0,T]$,
$$
   \|\vp\|_{U_T}^2 = \int_0^T \|\vp(s)\|_U^2 \; ds
=\sum_{j=1}^\infty \int_0^T \langle \vp(s), e_j\rangle_U^2\, ds.
$$
In particular, for any $j\geq 1$, the function $s \mapsto \langle
\vp(s),e_j\rangle_U$ belongs to $L^2([0,T];\re)$. Thus, it follows
from Step 2 that
$$\vp^n(\cdot):= \sum_{j=1}^n \langle \vp(\cdot),e_j\rangle_U \;
e_j$$ belongs to $C$. Moreover, it is straightforward to verify that
$\|\vp-\vp^n\|_{U_T}^2\to 0$ as $n\to \infty$.
This shows that $\vp\in C$. \qed

\vspace{0.3cm}

Therefore, taking into account the above lemma, $W(\vp)$ can be
defined for all $\ffi\in U_T$ following the standard method for
extending an isometry. This establishes the following property.

\begin{proposition}\label{spcyl}
For $t \geq 0$ and $\vp \in U$, set $W_t(\vp) = W(1_{[0,t]}(\cdot)\vp(\star ))$. Then the process $W=\{W_t(\vp),\, t \geq 0,\, \vp \in U\}$ is a cylindrical Wiener process as defined in Section \ref{cyl}, with $V$ there replaced by $U$ and $Q=\mbox{Id}_U$. In particular, for any $\vp \in U$, $\{W_t(\vp),\, t\geq 0\}$ is a Brownian motion with variance $t \Vert \vp\Vert_U$ and for all $s,t \geq 0$ and $\vp,\psi \in U$, $E(W_t(\vp)W_s(\psi)) = (s \wedge t) \langle \vp,\psi\rangle_U$.
\end{proposition}

   With this proposition, it becomes possible to use the stochastic integral defined in Section \ref{integral}.
   This defines the stochastic integral $g \cdot W$ for all $g \in L^2(\Omega\times[0,T]; U) \equiv L^2(\Omega; U_T)$.
   By definition of $U$, the complete orthonormal basis $(e_j)_j$ in the definition of $g\cdot W$ can be chosen
   such that $(e_j)_j \subset \cs(\red)$.

   Before discussing this further, we first relate the statement of Proposition \ref{spcyl} to Walsh's theory of stochastic integrals with respect to martingale measures. Let us recall that Walsh's theory of stochastic integration is
based on the concept of \emph{martingale measure}, which is a
stochastic process of the form $\{ M_t(A),\, \tf_t,\,  t\in [0,T],\, A\in
\cb_b(\red)\}$, where $\cb_b(\red)$ denotes the set of bounded Borel
sets of $\red$, and $(\tf_t)_t$ is a filtration satisfying the {\em usual conditions}. For the precise definition of a martingale measure, we refer to \cite[Chapter 2]{walsh}. Hence, in order to use Walsh's construction,  one has first to extend the generalized random field $\{W(\vp),\, \vp\in
\cinf_0(\re_+\times \red)\}$ to a martingale measure. More
precisely, using an approximation procedure similar to the one used in Lemma \ref{lem2.2}, one extends the definition of $W$ to
indicator functions of bounded Borel sets in $\re_+\times\red$ (for details see
\cite{DF} or \cite[p.13]{quer}). Then one sets
\beq
M_t(A)=W(1_{[0,t]}(\cdot) 1_A(\star )),\qquad t\in [0,T],\; A\in \cb_b(\red).
\label{mm}
\eeq

Moreover, if we let $(\tf_t)_t$ be the filtration generated by
$\{M_t(A),\, A\in \cb_b(\red)\}$ (completed and made right-continuous), then the process $\{M_t(A),\, \tf_t,\;
t\in [0,T],\, A\in \cb_b(\red)\}$ defines a worthy martingale measure in the sense of Walsh \cite{walsh}. Its {\em covariation measure} is determined by
$$\langle M(A),M(B)\rangle_t = t \int_{\red} \int_{\red} 1_A (x) f(x-y) 1_B(y) dx dy,$$
$t\in [0,T]$, $A,B\in \cb_b(\red)$, and its {\em dominating measure} coincides with the covariance measure (see \cite{DF}).

   One easily checks that, for $\ffi\in \cs(\red)$,
$$
   W_t(\vp) = \int_{\re_+} \int_{\red} 1_{[0,t]}(s)\, \vp(x)\, M(ds,dx),
$$
where the integral on the right-hand side is Walsh's stochastic integral.

\subsection{The real-valued stochastic integral for spatially homogeneous noise}
\label{homo}

The aim of this section is to exhibit the relationship between the stochastic integral constructed in Section \ref{integral} and the random field approach of Walsh \cite{walsh} and Dalang \cite{dalang}. Recall that the stochastic integral with respect to $M$ defined in \cite{walsh} only allows
function-valued integrands, and this theory was
extended in \cite{dalang} in order to cover more
general integrands, such as certain processes with values in the space of (Schwartz) distributions. We are going to show that these two integrals can be interpreted in the context of Section \ref{integral}.


   Recall that Walsh's stochastic integral $g\cdot M$ is defined when $g \in \cp_+$, where $\cp_+$ is the set of predictable processes $(\om,t,x) \mapsto g(t,x;\om)$ such that
$$
   \Vert g \Vert_+^2 := E\left( \int_0^T dt \int_{\red}dx \int_{\red}dy\,  |g(t,x)|\, f(x-y)\, |g(t,y)|   \right) <\infty.
$$
For $g \in \cp_+$, we can consider that $g \in L^2(\Omega;U_T)$ and set
\beq
   \Vert g \Vert_0^2 := E(\Vert g \Vert^2_{U_T})= E\left( \int_0^T dt \int_{\red} dx \int_{\red} dy\, g(t,x) f(x-y) g(t,y)   \right).
\label{norm0}
\eeq

In \cite{dalang}, Dalang considered the set $\cp_0$, which is the completion with respect to $\Vert \cdot \Vert_0$ of the subset $\ce_0$ of $\cp_+$ that consists of functions $g(s,x; \omega)$ such that $x \mapsto g(s,x; \omega) \in \cs(\red)$, for all $s$ and $\omega$, and he defined the stochastic integral $g\cdot M$ for all $g \in \cp_0$.

  Finally, in order to use the stochastic integral of Section \ref{integral}, let $(e_j)_j \subset \cs(\red)$ be a complete orthonormal basis of $U$, and consider the cylindrical Wiener process $\{W_t(\vp)\}$ defined in Proposition \ref{spcyl}. For any predictable process $g\in $ $L^{2}(\Omega \times [0,T];U)$, the
stochastic integral of $g$ with respect $W$ is
\beq
g\cdot W =\int_0^T g_s\, dW_s := \sum_{j=1}^\infty \int_0^T \langle g_s
,e_j\rangle_U\, dW_s(e_j), \label{2} \eeq
and the isometry property
is given by
\beq\label{isomrda}
E\left( \left( g\cdot W\right)^{2}\right) = E\left( \left(\int_0^T g_s\, dW_s\right)^2 \right)=  E\left( \int_{0}^{T}\left\|
g_s\right\| _{U}^{2}\, ds\right) .
\eeq
We note that the right-hand side of (\ref{2}) is essentially the definition of $W(\vp)$ in \cite{MS}. We also use the notation
$$
   \int_0^T \int_{\mathbb{R}^{d}} g(s,y)\, W(ds,dy)
$$
instead of $\int_0^T g_s\, dW_s$.

\begin{proposition} \label{w-d} (a) If $g \in \cp_+$, then $g \in L^2(\Omega \times [0,T]; U)$ and $g\cdot M = g \cdot W$, where the left-hand side is a Walsh integral and the right-hand side is defined as in \eqref{2}.

   (b) If $g \in \cp_0$, then $g \in L^2(\Omega \times [0,T]; U)$ and $g\cdot M = g \cdot W$, where the left-hand side is a Dalang integral and the right-hand side is defined as in \eqref{2}.

\end{proposition}

\noindent{\emph{Proof.} Let us prove part (a) in the statement. We first observe that if $g\in \cp_+$, then
\begin{align}
\|g\|^2_{L^2(\Om\times [0,T];U)} & = E\left( \int_0^T dt \int_{\red} dx \int_{\red} dy\, g(t,x) f(x-y) g(t,y)   \right) \nonumber \\
& \leq \|g\|^2_+ < +\infty.
\label{wd}
\end{align}
This implies that $g\in L^2(\Om\times [0,T];U)$.

   Secondly, in order to check the equality of the integrals, we use the fact that the set of elementary processes is dense in $(\cp_+,\|\cdot\|_+)$ (see \cite[Proposition 2.3]{walsh}). Hence, by inequality (\ref{wd}), it suffices to show that both
integrals coincide when $g$ is an elementary process of the form
\beq g(t,x;\om)= 1_{(a,b]}(s) 1_A(x) X(\om), \label{6} \eeq where
$0\leq a<b\leq T$, $A\in \cb_b(\red)$ and $X$ is a bounded and
$\tf_a$-measurable random variable.

On one hand, when $g$ has the particular form \eqref{6}, according to \cite{walsh} and \eqref{mm},
\begin{align*}
\int_0^T \int_{\red} g(t,x)\, M(dt,dx) & = \left[ M_b(A)-M_a(A)\right] X \\
& = \left[ W(1_{(0,b]}(\cdot) 1_A(\star ))-W(1_{(0,a]}(\cdot) 1_A(\star ))\right] X \\
& = W(1_{(a,b]}(\cdot) 1_A(\star )) X.
\end{align*}
On the other hand, by the very definition of the integral (\ref{2}),
\begin{align*}
\int_0^T g_t\, dW_t & = \sum_{j=1}^\infty \int_a^b X \langle 1_A,e_j\rangle_{U}\, dW_t(e_j)\\
& = X \sum_{j=1}^\infty \langle 1_A,e_j\rangle_{U} \left[ W_b(e_j)-W_a(e_j)\right]\\
& = X \sum_{j=1}^\infty \langle 1_A,e_j\rangle_{U}\, W(1_{(a,b]}(\cdot)e_j)\\
& = X W(1_{(a,b]}(\cdot)1_A(\star )),
\end{align*}
which implies that
$$ \int_0^T \int_{\red} g(t,x) M(dt,dx) = \int_0^T g_t\, dW_t,$$
for all $g$ of the form (\ref{6}). This concludes the first part of the proof.

   Concerning part (b), let us point out that $\cp_0$ is the completion of $\ce_0$ with respect to $\|\cdot\|_0$ (see (\ref{norm0})), where the latter coincides with the norm in $L^2(\Om\times [0,T];U)$ for smooth elements. Hence, since $\ce_0\subset \cp_+\subset L^2(\Om\times [0,T];U)$, any $\|\cdot\|_0$-limit $g$ of a sequence $(g_n)_n\subset \ce_0$ will determine a well-defined element in  $L^2(\Om\times [0,T];U)$.

Moreover, as a consequence of this, we will only need to check the equality of the integrals for integrands $g$ in $\ce_0$. Since such elements are contained in $\cp_+$, Dalang's integral of $g$ with respect to the martingale measure $M$ turns out to be a Walsh integral, so that we can conclude by using the first part of the proof.
\qed

\begin{remark} According to Proposition \ref{w-d}, when one integrates an element of $\cp_+$, it is possible to use either the Walsh integral or the integral with respect to a cylindrical Wiener process. However, the Walsh integral enjoys additional properties, in part because it is possible to make use of the dominating measure, which can be very useful in certain estimates.  For example, establishing H\"older continuity of the solution to the $1$-dimensional stochastic wave equation, in which a Walsh integral appears, is an easy exercise \cite[Exercise 3.7]{walsh}, while for the $3$-dimensional stochastic wave equation, this is quite involved \cite{ds2}.
\end{remark}

\subsection{Examples of integrands}
\label{distri}

In this section, we aim to provide useful examples of random
distributions which belong to $L^2(\Om\times [0,T];U)$, that is, for which we can define the stochastic integral (\ref{2}) with respect to $W$.

   Recall that an element $\Theta\in \cs'(\red)$ is a
\emph{non-negative distribution with rapid decrease} if $\Theta$ is a non-negative measure and if
$$\int_{\red} (1+|x|^2)^{k/2}\, \Theta(dx) < +\infty,$$
for all $k>0$ (see \cite{schwartz}).

   Recall that $\mu$ is the spectral measure of $W$. We consider the
following hypothesis.

\begin{hypothesis}\label{hypA}
Let $\Gam$ be a function defined on $\re_+$ with values in $\cs'(\red)$ such that, for all
$t>0$, $\Gamma(t)$ is a non-negative distribution with rapid
decrease, and
\beq \int_0^T dt \int_{\mathbb{R}^{d}} \mu(d\xi)\, |\mathcal{F}
\Gamma(t)(\xi)|^2 < \infty. \label{mu-L}
\eeq
In addition, $\Gamma$ is a non-negative measure of the form $\Gamma(t, dx)dt$
such that, for all $T>0$,
$$
\sup_{0\leq t\leq T} \Gamma(t,\mathbb{R}^{d}) <\infty.
$$
\end{hypothesis}

\vspace{0.3cm}

The main examples of integrands are provided by the following
proposition (see \cite[Proposition 3.3 and Remark 3.4]{nualartquer}). In comparison with the analogous result by Dalang \cite[Theorem 2]{dalang}, Proposition \ref{prop1} does not require that the stochastic process $Z$ have a spatially homogeneous covariance (see Hypothesis A in \cite{dalang}).

\begin{proposition}\label{prop1}
Assume that $\Gamma$ satisfies Hypothesis \ref{hypA}. Let
$Z=\{Z(t,x),\, (t,x)\in [0,T]\times \mathbb{R}^{d}\}$ be a predictable
process such that
\beq \sup_{(t,x)\in \lbrack 0,T]\times
\mathbb{R}^{d}}E(|Z(t,x)|^{p})< \infty, \label{sup}
\eeq
for some
$p\geq 2$. Then, the random measure $G=\{G(t,dx)=Z(t,x)\Gamma
(t,dx),\; t\in [0,T]\}$ is a predictable process with values in $L^{p}(\Omega \times [0,T];U)$. Moreover, 
$$
   E\left( \|G\|^2_{U_T}\right) = E\left[  \int_0^T dt \int_{\red} \mu(d\xi)\, |\tf (\Gam(t)Z(t))(\xi)|^2\right]
$$
and
$$
E\left( \|G\|^p_{U_T}\right) 
 \leq  C \int_0^T  dt\left( \sup_{x\in \red} E(|Z(t,x)|^p) \right)   \int_{\mathbb{R}^d}\mu(d\xi)\, |\mathcal{F} \Gamma(t) (\xi)|^2   .
$$
\end{proposition}

   The integral of $G=\{G(t,dx)=Z(t,x)\Gamma
(t,dx),\; t\in [0,T]\}$ with respect to $W$ will be also denoted by
\beq
G\punt W =\int_0^T \int_{\red} \Gam(s,y)Z(s,y)W(ds,dy).
\label{int-gz}
\eeq

It is worth pointing out two key steps in the proof of this
proposition (see \cite{nualartquer}): the first is to check that under Hypothesis \ref{hypA}, $\Gam$ belongs to $U_T=L^2([0,T];U)$; the second is to notice that if $\Gam$ and $Z$ satisfy, respectively, Hypothesis \ref{hypA} and condition (\ref{sup}), then $G(t)=Z(t,\star )\Gam(t,\star )$ defines a distribution with rapid decrease, almost surely.

\begin{remark}\label{conus-dalang}
 We note that \cite{CD} presents a further extension of Walsh's stochastic integral, with which it becomes possible to integrate certain random elements of the form $Z(t,z)\Gam(t,\star)$, where $\Gam$ is a tempered distribution which is not necessarily non-negative. This extension is useful for studying the stochastic wave equation in high spatial dimensions. 
\end{remark}

\subsection{The Dalang-Mueller extension of the stochastic integral}
\label{ext-dm}

We briefly summarize here the function-valued stochastic integral constructed in \cite{dm}. This is an extension of Walsh's stochastic integral, where one
integrates processes that take values in $L^2(\red)$ (or a weighted $L^2$-space) and the
value of the integral is in the same $L^2$-space.

Suppose that $s\mapsto \Gam(s)\in \cs'(\red)$ satisfies:
\begin{itemize}
    \item[1)] For all $s\geq 0$, $\tf \Gam(s)$ is a function and
    $$\int_0^T ds \sup_{\xi\in \red} \int_{\red} \mu(d\eta)\, |\tf \Gam(s)(\xi-\eta)|^2 <+\infty.$$
    \item[2)] For all $\phi\in \cinf_0(\red)$, $\sup_{0\leq s\leq T} \Gam(s) *\phi$ is a bounded function on $\red$.
\end{itemize}
Suppose that $s\mapsto Z(s)\in L^2(\red)$ satisfies:
\begin{itemize}
    \item[3)] For $0\leq s\leq T$, $Z(s)\in L^2(\red)$ a.s., $Z(s)$ is $\tf_s$-measurable, and $s\mapsto Z(s)$ is mean-square continuous from $[0,T]$ into $L^2(\red)$.
\end{itemize}
For such $\Gam$ and $Z$, one sets
\beq
   I_{\Gam,Z}:=\int_0^T ds \int_{\red} d\xi \, E\left( |\tf Z(s)(\xi)|^2\right) \int_{\red} \mu(d\eta)\,
   |\tf \Gam(s)(\xi-\eta)|^2 <+\infty.
\label{iz1} \eeq
Then the stochastic integral
\beq
v_{\Gam,Z}=\int_0^T \int_{\red} \Gam(s,\star -y) Z(s,y)\, M(ds,dy)
\label{102}
\eeq
is defined as an element of $L^2(\Om\times \red,dP\times dx)$, such that
\beq
E\left(\|v_{\Gam,Z}\|^2_{L^2(\red)} \right) = I_{\Gam,Z}.
\label{112}
\eeq
This definition is obtained in three steps.\\

   {a)} If, in addition to 1), $\Gam(s)\in \cinf(\red)$, for $0\leq s\leq T$, and in
addition to 3), $Z(s)\in \cinf_0(\red)$ and there is a compact $K\subset \red$
    such that $\mbox{supp } Z(s)\subset K$, for $0\leq s\leq T$, then
$$v_{\Gam,Z}(x)=\int_0^T \int_{\red} \Gam(s,x-y) Z(s,y)M(ds,dy),$$
where the right-hand side is a Walsh stochastic integral. Equality (\ref{112}) is checked
by direct calculation (see \cite[Lemma 1]{dm}).

   {b)} If $\Gam$ is as in {a)} and $Z$ satisfies 3), then one checks that
$$
   \lim_{m\rightarrow \infty}\lim_{n\rightarrow \infty}
I_{\Gam,Z-(Z1_{[-m,m]})*\psi_n}=0,
$$
where $(\psi_n) \subset C^\infty_0(\red)$ is a sequence that converges to the Dirac distribution,
 and one sets
$$v_{\Gam,Z}= \lim_{m\rightarrow \infty}\lim_{n\rightarrow \infty}
v_{\Gam,(Z1_{[-m,m]})*\psi_n},$$ where the limits are in $L^2(\Om\times \red,dP\times
dx)$.

   {c)} If $\Gam$ satisfies 1) and 2), and $Z$ satisfies 3), then one checks that
$$\lim_{n\rightarrow \infty} I_{\Gam-\Gam*\psi_n,Z}=0$$
and one sets
$$v_{\Gam,Z}=\lim_{n\rightarrow \infty} v_{\Gam*\psi_n,Z},$$
where the limit is in $L^2(\Om\times \red,dP\times dx)$: see \cite[Theorem 6]{dm}.\\

In comparison with the stochastic integral 
of Section \ref{homo}, we remark that the process $Z$ verifies $\sup_{s\in
[0,T]} E(\|Z(s)\|^2_{L^2(\red)}) <+\infty$, rather than
\eqref{sup}, and the resulting integral $v_{\Gam,Z}$, as a random function of $x$, belongs to $L^2(\Om\times\red)$.

  We now relate this stochastic integral to the one defined in Section \ref{homo}.

\begin{proposition}\label{integrals}
Assume that $\Gam$ and $Z$ satisfy conditions 1), 2) and 3) above. Then:
\begin{itemize}
 \item[(i)] For almost all $x\in \red$, the element $\Gam(\cdot,x-\star ) Z(\cdot,\star )$ belongs to $L^2(\Om\times [0,T];U)$. Hence, as in \eqref{2}, we can define the (real-valued) stochastic integral
$$\ci_{\Gam,Z}(T,x):=\int_0^T \int_{\red} \Gam(s,x-y) Z(s,y)\, W(ds,dy), \qquad \text{for a.a.}\; x\in\red.$$

\item[(ii)] $\ci_{\Gam,Z}(T,\star )\in L^2(\Om\times \red)$ and
$\|\ci_{\Gam,Z}(T,\star ) \|^2_{L^2(\Om\times \red)} = I_{\Gam,Z}$.

\item[(iii)] $\ci_{\Gam,Z}(T,\star )=v_{\Gam,Z}$ in $L^2(\Om\times \red)$.
\end{itemize}
\end{proposition}

\noindent \emph{Proof.} We will split the proof in three steps, which essentially correspond to  the construction
of the Dalang-Mueller integral $v_{\Gam,Z}$.

\medskip

\noindent {\it{Step 1.}} Let us assume first that $\Gam$ and $Z$ satisfies the hypotheses in a) above. Then, as we pointed out there, for all $x\in \red$, the stochastic integral $v_{\Gam,Z}(x)$ can be defined
as a Walsh stochastic integral. Hence, by Proposition \ref{w-d}(a), the integrand $(s,y)\mapsto \Gam(s,x-y) Z(s,y)$ defines
an element in $L^2(\Om\times [0,T]; U)$ and, for all $x\in \red$, $v_{\Gam,Z}(x)=\ci_{\Gam,Z}(T,x)$. Condition (ii) in
the statement can be deduced from this latter equality and \eqref{112}.

\medskip

\noindent {\it{Step 2.}} Assume now that $\Gam$ is as in Step 1 and $Z$ satisfies condition 3). Then, as in b) above,
there exists a sequence of processes $(Z_n)_n$ such that, for all $n\geq 1$, $Z_n$
satisfies the hypotheses in a) and $I_{\Gam,Z_n-Z}$ converges to zero as $n$ tends to infinity. For this sequence,
\beq\label{104}
    v_{\Gam,Z}:= \lim_{n\rightarrow \infty} v_{\Gam,Z_n}
       = \lim_{n\rightarrow \infty} \ci_{\Gam,Z_n}(T,\star )
\eeq
by Step 1,
where the limit is in $L^2(\Om\times \red)$.

   We now check property (i) in the statement of the proposition. Observe that, by Proposition \ref{prop1},
\begin{align}\label{104a}
& \int_{\red} dx\, \|\Gam(\cdot,x-\star )[Z_n(\cdot,\star )-Z(\cdot,\star )]\|^2_{L^2(\Om\times [0,T];U)} \\
& \qquad = \int_{\red} dx\,E\left( \int_0^T ds \int_{\red} \mu(d\eta)\, |\tf\big(\Gam(s,x-\star )[Z_n(s,\star )-Z(s,\star )]\big)(\eta)|^2 \right).
\nonumber
\end{align}
Use the very last lines in the proof of \cite[Lemma 1]{dm} to see that this is equal to $I_{\Gam,Z_n-Z}$.
Since this quantity converges to zero as $n\to \infty$, we deduce that there exists a subsequence 
$(n_j)_j$ such that, for almost all $x\in \red$, 
$$
   \lim_{j\rightarrow \infty} \left\| \Gam(\cdot,x-\star )Z_{n_j}(\cdot,\star )-\Gam(\cdot,x-\star )Z(\cdot,\star )\right\|_{L^2(\Om\times [0,T]; U)}=0.
$$
This implies that, for almost all $x\in \red$, the element $(s,y)\mapsto \Gam(s,x-y)Z(s,y)$ belongs to $L^2(\Om\times [0,T]; U)$, and we can define the (real-valued) stochastic integral 
\beq
\ci_{\Gam,Z}(T,x):=\int_0^T \int_{\red} \Gam(s,x-y) Z(s,y)\, W(ds,dy),
\label{103}
\eeq
and
$$
   \ci_{\Gam,Z}(T,x) = \lim_{j \to \infty} \ci_{\Gam,Z_{n_j}}(T,x) \qquad\mbox{in } L^2(\Omega).
$$

   Notice that 
\beq\label{103a}
   \| \ci_{\Gam,Z_n}(T,\star ) - \ci_{\Gam,Z}(T,\star ) \|_{L^2(\Om\times \red)}^2 = \|\ci_{\Gam,Z_n - Z}(T,\star ) \|_{L^2(\Om\times \red)}^2.
\eeq
By the isometry property \eqref{isomrda}, this is equal to \eqref{104a}, and therefore to $I_{\Gamma, Z_n - Z}$, which tends to $0$ as $n \to \infty$. Therefore, using Step 1, we see that
$$
   \| \ci_{\Gam,Z}(T,\star ) \|_{L^2(\Om\times \red)}^2 = \lim_{n\to\infty} \| \ci_{\Gam,Z_n}(T,\star )\|_{L^2(\Om\times \red)}^2 = \lim_{n\to\infty} I_{\Gamma,Z_n} = I_{\Gamma,Z},
$$
which proves (ii). The arguments following \eqref{103a} and \eqref{104a} prove (iii).
\medskip

\noindent {\it{Step 3.}} In this final part, we assume that $\Gam$ and $Z$ satisfy conditions 1), 2) and 3). 
Then, it is a consequence of step c) above 
that there exists $(\Gam_n)_n$ such that, for all $n\geq 1$, $\Gam_n$ verifies the assumptions 
of the previous step and 
$$\lim_{n\rightarrow \infty} I_{\Gam_n-\Gam,Z}=0.$$
In order to prove parts (i), (ii) and (iii) for this case, one can follow exactly the same lines as we have done in Step 2. We omit the details.
\qed
\vskip 16pt

As we will explain in Section \ref{relation-dm}, for the particular case of the stochastic wave equation, it is useful to consider
stochastic integrals of the form $v_{\Gam,Z}$ which take values in some weighted $L^2$-space. We now describe this situation.

Fix $k>d$ and let $\theta :\red \rightarrow \re$ be a smooth function for which there are constants
$0<c<C$ such that
$$c(1 \wedge |x|^{-k}) \leq \theta(x) \leq C (1 \wedge |x|^{-k}).$$
The weighted $L^2$-space $L^2_\theta$ is the set of measurable $g:\red \rightarrow \re$ such that $\|g\|_\theta <+\infty$, where
$$\|g\|_\theta^2 = \int_{\red} |g(x)|^2\, \theta(x)\, dx.$$
Consider a function $s\mapsto \Gam(s)\in \cs'(\red)$ that satisfies 1), 2) above, and, in
addition,
\begin{itemize}
    \item[4)] There is $R>0$ such that for $s\in [0,T]$, $\mbox{supp } \Gam(s) \subset B(0,R)$.
\end{itemize}
For a stochastic process $Z$, we consider the following hypothesis:
\begin{itemize}
    \item[5)] For $0\leq s\leq T$, $Z(s)\in L^2_\theta$ a.s., $Z(s)$ is $\tf_s$-measurable, and $s\mapsto Z(s)$ is mean-square continuous from $[0,T]$ into $L^2_\theta$.
\end{itemize}
Then the stochastic integral
\beq\label{103b}
v^\theta_{\Gam,Z}=\int_0^T \int_{\red} \Gam(s,\star -y) Z(s,y)M(ds,dy)
\eeq
is defined as an element of $L^2(\Om\times \red,dP\times \theta(x)dx)$, such that
$$E(\|v^\theta_{\Gam,Z}\|^2_{L^2_\theta}) \leq I^\theta_{\Gam,Z},$$
where
$$I^\theta_{\Gam,Z}:=\int_0^T ds \, E(\|Z(s,\star )\|^2_{L^2_\theta}) \sup_{\xi\in \red}  \int_{\red}
\mu(d\eta)\,  |\tf \Gam(s)(\xi-\eta)|^2.$$
This definition is obtained by showing
that $Z_n(s,\star ):= Z(s,\star ) 1_{[-n,n]} (\star )$ also satisfies 5) as well as 3). Therefore,
$v^\theta_{\Gam,Z_n}=v_{\Gam,Z_n}$ is defined as an element of $L^2(\Om\times
\red,dP\times dx)$, and one checks that this element also belongs to $L^2(\Om\times
\red,dP\times \theta(x)dx)$, and
$$\lim_{n\rightarrow \infty} I^\theta_{\Gam,Z-Z_n}=0,$$
provided that $\Gam$ satisfies 1), 2) and 4). Then one sets
$$v^\theta_{\Gam,Z} = \lim_{n\rightarrow \infty} v_{\Gam,Z_n},$$
where the limit is in $L^2(\Om\times \red,dP\times \theta(x)dx)$: see \cite[Theorem 12]{dm}.


\section{Infinite-dimensional integration theory}
\label{infdim}

In this section, we first sketch the construction of the
infinite dimensional stochastic integral in the setup of Da Prato and
Zabczyk in \cite{dz}. For this, we will define the general concept
of Hilbert-space-valued $Q$-Wiener process and study its relationship
with the cylindrical Wiener process considered in Section
\ref{integral}. Then we will show that the stochastic integral
constructed in Section \ref{integral} can be inserted into this more
abstract setting. In particular, we will treat specifically the case of the
standard cylindrical Wiener process given by the spatially
homogeneous noise
described in Section \ref{prel}.

   We begin by recalling some facts concerning nuclear and
Hilbert-Schmidt operators on Hilbert spaces.

\subsection{Nuclear and Hilbert-Schmidt operators}
\label{recall}

Let $E,G$ be Banach spaces and let $L(E,G)$ be the vector space of all linear bounded operators from $E$ into $G$.
We denote by $E^*$ and $G^*$ the dual spaces of $E$ and $G$, respectively.

   An element $T\in L(E,G)$ is said to be a \emph{nuclear} operator if
there exist two sequences $(a_j)_j\subset G$ and $(\vp_j)_j\subset
E^*$ such that
$$
T(x) =\sum_{j=1}^\infty a_j\, \vp_j(x), \qquad\mbox{for all } x\in E,
$$
and
$$
\sum_{j=1}^\infty \|a_j\|_G\, \|\vp_j\|_{E^*} < +\infty.
$$

The space of all nuclear operators from $E$ into $G$ is denoted by
$L_1(E,G)$. When endowed with the norm
$$
\|T\|_1 = \inf \left\{ \sum_{j=1}^\infty \|a_j\|_G\,
\|\vp_j\|_{E^*}:\;  T(x) =\sum_{j=1}^\infty a_j \vp_(x), \; x\in E
\right\},
$$
it is a Banach space.

   Let $H$ be a separable Hilbert space and let
$(e_j)_j$ be a complete orthonormal basis in $H$. For $T\in
L_1(H,H)$, the {\em trace} of $T$ is
\beq
  \mbox{Tr}\; T = \sum_{j=1}^\infty \langle T( e_j), e_j\rangle_H.\label{trace}
\eeq
One proves that if $T\in L_1(H):=L_1(H,H)$, then $\mbox{Tr}\; T$ is a
well-defined real number and its value does not depend on the
choice of the orthonormal basis (see, for instance,
\cite[Proposition C.1]{dz}). Further, according to \cite[Proposition
C.3]{dz}, a non-negative definite operator $T\in L(H)$ is nuclear if and only
if, for an orthonormal basis $(e_j)_j$ on $H$,
$$
\sum_{j=1}^\infty \langle T( e_j),e_j\rangle_H < +\infty.
$$
Moreover, in this case, $\mbox{Tr}\; T= \|T\|_1$.

Let $V$ and $H$ be two separable Hilbert spaces and $(e_k)_k$ a
complete orthonormal basis of $V$. A bounded linear operator
$T:V\rightarrow H$ is said to be \emph{Hilbert-Schmidt} if
$$\sum_{k=1}^\infty \|T(e_k)\|^2_H < +\infty.$$
It turns out that the above property is independent of the choice of
the basis in $V$. The set of Hilbert-Schmidt operators from $V$ into
$H$ is denoted by $L_2(V,H)$. The norm in this space is defined by
\beq \|T\|_2 = \left( \sum_{k=1}^\infty \|T(e_k)\|^2_H
\right)^{1/2},\label{norm-hs} \eeq and defines a Hilbert space with
inner product
\beq \langle S, T\rangle_2 = \sum_{k=1}^\infty \langle
S(e_k),T(e_k)\rangle_H. \label{escalar}
\eeq
Finally, let us point out that (\ref{trace}) and (\ref{norm-hs}) imply that if $T\in L_2(V,H)$, then $TT^* \in L_1(H)$, where $T^*$ is the adjoint operator of $T$, and
\beq
   \|T\|_2^2 = \mbox{Tr}\; (T T^*).
\label{adjoint}
\eeq

We conclude this section by recalling the definition and some
properties of the pseudo-inverse of bounded linear operators (see, for instance, \cite[Appendix C]{pr}).

Let $T\in L(V,H)$ and $\mbox{Ker }T := \{x\in V: T(x)=0\}$. The
\emph{pseudo-inverse} of the operator $T$ is defined by
$$
   T^{-1}:= \left( T_{|_{(\mbox{\scriptsize Ker }T)^\bot}}\right)^{-1} :\ T(V)\to (\mbox{Ker }T)^\bot.
$$
Notice that $T$ is one-to-one on $(\mbox{Ker }T)^\bot$ (the orthogonal complement of $\mbox{Ker }T$) and $T^{-1}$ is linear and bijective.

If $T\in L(V)$ is a bounded linear operator defined on $V$ and
$T^{-1}$ denotes the pseudo-inverse of $T$, then (see \cite[Proposition
C.0.3]{pr}):
\begin{enumerate}
  \item $( T(V), \langle\punt ,\punt\rangle_{T(V)} )$ defines a Hilbert
  space, where
  $$\langle x,y\rangle_{T(V)} := \langle
  T^{-1}(x),T^{-1}(y)\rangle_V,\qquad x,y\in T(V).$$
  \item Let $(e_k)_k$ be an orthonormal basis of
  $(\mbox{Ker }T)^\bot$. Then $(T(e_k))_k$ is an orthonormal basis
  of $( T(V), \langle\punt ,\punt \rangle_{T(V)} )$.
\end{enumerate}
Finally, according to \cite[Corollary C.0.6]{pr}, if $T\in L(V,H)$
and we set $Q:= T T^* \in L(H)$, then we have $\mbox{Im } Q^{1/2} =
\mbox{Im } T$ and
$$
\left\| Q^{-1/2}(x)\right\|_H = \|T^{-1} (x)\|_V, \qquad x\in \mbox{Im
} T,
$$
where $Q^{-1/2}$ is the pseudo-inverse of $Q^{1/2}$.

\subsection{Hilbert-space-valued Wiener processes}
\label{wiener}

   The stochastic integral presented in Da Prato and Zabczyk \cite{dz} is defined with respect to a class of Hilbert-space-valued processes, namely $Q$-Wiener processes, which we now introduce.

We consider a separable Hilbert space $V$ and a linear, symmetric (self-adjoint)
non-negative definite and bounded operator $Q$ on $V$ such that $\mbox{Tr}\; Q <+\infty$.
\begin{definition}\label{Qwiener}
A $V$-valued stochastic process $\{\cw_t,\ t\geq 0\}$ is called a
{\itshape{$Q$-Wiener process}} if
(1) $\cw_0=0$,
(2) $\cw$ has continuous trajectories,
(3) $\cw$ has independent increments,
and (4) the law of $\cw_t-\cw_s$ is Gaussian with mean zero and covariance operator $(t-s)Q$, for all $0\leq s\leq t$.
\end{definition}
We recall that according to \cite[Section 2.3.2]{dz}, condition (4)
above means that for any $h\in V$ and $0\leq s\leq t$, the
real-valued random variable $\langle \cw_t-\cw_s, h\rangle_V$ is
Gaussian, with mean zero and variance $(t-s) \langle
Qh,h\rangle_V$. In particular, using \eqref{trace}, we see that $E(\Vert \cw_t\Vert_V^2) = t\, \mbox{Tr}\;  Q$, which is one reason why the assumption $\mbox{Tr}\;  Q < \infty$ is essential.

   Let $(e_j)_j$ be an orthonormal basis of $V$ that consists of eigenvectors of $Q$ with corresponding eigenvalues $\lambda_j$, $j \in \N^*$.
   Let $(\beta_j)_j$ be a sequence of independent real-valued standard Brownian motions on a probability space $(\Omega,\cf,P)$. Then the $V$-valued process
\begin{equation}\label{wienerrepr}
   \cw_t = \sum_{j=1}^\infty \sqrt{\lambda_j}\, \beta_j(t) e_j
\end{equation}
(where the series converges in $L^2(\Omega;\cac([0,T];V))$), defines a $Q$-Wiener process
on $V$ (see (2.1.2) in \cite{pr}). We note that $\sqrt{\lambda_j}\, e_j = Q^{1/2} (e_j)$.
In the special case where $V$ is finite-dimensional, say dim~$V = n$, then $Q$ can be
identified with an $n\times n$-matrix  which is the variance-covariance matrix of
$\{\cw_t\}$, and $\{\cw_t\}$ has the same law as $\{Q^{1/2} \cw^0_t \}$, where $\{\cw^0_t
\}$  is a standard Brownian motion with values in $\R^n$.

If $\{\cw_t,\ t\geq 0\}$ is a $Q$-Wiener process on $V$, there is a
natural way to associate to it a cylindrical Wiener process in the
sense of Definition \ref{cyl}. Namely, for any $h\in V$ and $t\geq
0$, we set $W_t(h):=\langle \cw_t , h\rangle_V$. Using polarization, one checks that $\{W_t(h),\
t\geq 0,\ h\in V\}$ is a cylindrical Wiener process on $V$ with
covariance operator $Q$.  Note that in this case, $W_t(e_j) = \sqrt{\lambda_j}\, \beta_j(t)$, so the Brownian motions $\beta_j$ in \eqref{wienerrepr} are given by $\beta_j(t) = W_t(v_j)$, where
\beq\label{e3.6}
   v_j = \lambda_j^{-1/2} e_j = Q^{-1/2}(e_j),\qquad \mbox{for } j \geq 1 \mbox{ with } \lambda_j \neq 0.
\eeq
In particular, $(v_j)_j$ is a complete orthonormal basis of the space $V_Q$ of Section \ref{integral}.

   However, it is not true in general that any cylindrical Wiener process is associated to a $Q$-Wiener process on a Hilbert space. Indeed, we have the following result (see
\cite[p.177]{metivier}).

\begin{theorem}\label{t1}
Let $V$ be a separable Hilbert space and $W$ a cylindrical Wiener
process on $V$ with covariance $Q$. Then, the following three
conditions are equivalent:
\begin{enumerate}
  \item $W$ is associated to a $V$-valued $Q$-Wiener process $\cw$, in the sense
that $\langle \cw_t,h\rangle_V = W_t(h)$, for all $h\in V$.
  \item For any $t\geq 0$, $h \mapsto W_t(h)$ defines a Hilbert-Schmidt operator from $V$ into
$L^2(\Om, \tf,P)$.
  \item  $\mbox{\rm Tr}\; Q <+\infty$.
\end{enumerate}
If any one of the above conditions holds, then the norm of the Hilbert-Schmidt operator $h\mapsto W_t(h)$, as an element of $L_2\left( V,L^2(\Om, \tf,P)\right)$, is given by
$$\|W_t\|_2 =E(\|\cw_t\|^2_V)= t \mbox{ \rm Tr}\; Q.$$
\end{theorem}

As a consequence of the above result, if $\mbox{dim } V = +\infty$ and if $W$ is a
standard cylindrical Wiener process on $V$, that is $Q=\mbox{Id}_V$, then there is no
$Q$-Wiener process $\cw$ associated to $W$. However, as we will explain in Section
\ref{cylwin}, it will be possible to find a Hilbert-space-valued Wiener process  with
values in a larger Hilbert space $V_1$ which will correspond to $W$ in a certain sense.

\subsection{$H$-valued stochastic integrals}
\label{daprato}

We now sketch the construction of the infinite-dimensional
stochastic integral of \cite{dz}. Let $V$ and $H$ be two separable Hilbert spaces and let $\{\cw_t,\ t\geq 0\}$ be a $Q$-Wiener process defined on $V$. We note by $(\tf_t)_t$ the (completed) filtration generated by $\cw$. In \cite{dz}, the objective is to construct the
$H$-valued stochastic integral
$$
   \int_0^t \Phi_s\, d\cw_s,\qquad t\in [0,T],
$$
where $\Phi$ is a process with values in the space of linear but not necessarily bounded operators from $V$ into $H$.

Consider the subspace $V_0:=Q^{1/2}(V)$ of $V$ which, endowed with
the inner product
$$\langle h,g\rangle_0 := \langle Q^{-1/2} h, Q^{-1/2} g\rangle_V,$$
is a Hilbert space. Here $Q^{-1/2}$ denotes the pseudo-inverse of
the operator $Q^{1/2}$ (see Section \ref{recall}). Let us also
set
$$
   L_2^0 := L_2(V_0,H),
$$
which is the Hilbert space of all Hilbert-Schmidt operators from $V_0$ into $H$, equipped, as in \eqref{escalar}, with the inner product
\beq\label{L02norm}
\langle \Phi, \Psi\rangle_{L_2^0}= \sum_{j=1}^\infty \langle \Phi
\tilde e_j, \Psi \tilde e_j\rangle_H,\qquad \Phi, \Psi\in L_2^0,
\eeq
where $(\tilde e_j)_j$ is any complete orthonormal basis of $V_0$. In particular, using the fact that we can take
\beq\label{e3.8}
   \tilde e_j = \sqrt{\lam_j}\, e_j = Q^{1/2} (e_j), \qquad j\geq 1,\ \lam_j>0,
\eeq
where the $(e_j)_j$ are as in \eqref{wienerrepr} (see condition 2.~in the final part of Section \ref{recall}) and
applying (\ref{adjoint}), the norm of $\Psi \in L_2^0$ can be expressed as
$$\|\Psi\|^2_{L_2^0}= \| \Psi\circ Q^{1/2} \|_{L_2(V,H)}^2 = \mbox{Tr }(\Psi Q \Psi^*).$$
We note that in the case where dim~$V = n < +\infty$ and dim~$H=m < +\infty$, then it is natural to identify $\Psi \in L_2^0$ with an $m\times n$-matrix and $Q$ with an $n\times n$-matrix. The norm of $\Psi$ corresponds to a classical matrix norm of $\Psi Q^{1/2}$ (whose square is the sum of squares of entries of $\Psi Q^{1/2}$).

Let $\Phi=\{\Phi_t,\ t\in [0,T]\}$ be a measurable $L_2^0$-valued process. We define the norm of $\Phi$ by
$$
   \| \Phi\|_T :=\left[ E\left(\int_0^T \|\Phi_s\|^2_{L_2^0} \; ds \right)\right]^{1/2}.
$$
The aim of \cite[Chapter 4]{dz}, is to define the stochastic integral with respect to
$\cw$ of any $L_2^0$-valued predictable process $\Phi$ such that $\|\Phi\|_T<\infty$.
More precisely, Da Prato and Zabczyk first consider \textit{simple} processes, which are
of the form $\Phi_t=\Phi_0 1_{(a,b]}(t)$, where $\Phi_0$ is any $\tf_a$-measurable
$L(V,H)$-valued random variable and $0\leq a<b\leq T$. For such processes, the stochastic
integral takes values in $H$ and is defined by the formula
\begin{equation}\label{simple_int}
\int_0^t \Phi_s\, d\cw_s := \Phi_0 (\cw_{b\wedge t} - \cw_{a\wedge t}), \qquad t\in [0,T].
\end{equation}
The map $\Phi \mapsto \int_0^\cdot \Phi_s d\cw_s$ is an isometry
between the set of simple processes and the space $\mathcal{M}_H$ of
square-integrable $H$-valued $(\tf_t)$-martingales $X=\{X_t,\; t\in
[0,T]\}$ endowed with the norm $\|X\|=[ E (\|X_T\|^2_H) ]^{1/2}$.
Indeed, as it is proved in \cite{dz} (see also \cite[Proposition
2.3.5]{pr}), the isometry property for simple processes reads \beq
E\left( \left\| \int_0^T \Phi_t\, d\cw_t \right\|_H^2 \right) = \|
\Phi\|_T^2.\label{iso} \eeq

\begin{remark} The appearance of $\|\cdot\|_T$ can be understood 
by considering the case where $\Phi(t)=\Phi_0 1_{(a,b]}(t)$, where $\Phi_0\in L(V,H)$ is deterministic and $0\leq a<b\leq T$. Indeed, in this case, using \eqref{simple_int} and the representation \eqref{wienerrepr},
$$
   E\left( \left\| \int_0^T \Phi_t\, d\cw_t \right\|_H^2 \right) = E\left( \left\| \sum_j \sqrt{\lambda_j}\, (\beta_j(b) - \beta_j(a))\, \Phi_0(e_j)  \right\|_H^2 \right),
$$
and the right-hand side is equal to
\begin{align*}
   \sum_j \lambda_j\, (b-a) \left\|\Phi_0(e_j)  \right\|_H^2
   &= (b-a) \sum_j  \left\|\Phi_0(Q^{1/2}e_j)  \right\|_H^2
   = (b-a) \left\|\Phi_0 \circ Q^{1/2}\right\|_{L_2(V,H)}^2\\
   &= E\left(\int_0^T \left\|\Phi_s\right\|^2_{L_2^0} \; ds\right).
\end{align*}
\end{remark}

Once the isometry property \eqref{iso} is established, a completion argument is used to
extend the above definition to all $L_2^0$-valued predictable
processes $\Phi$ satisfying $\|\Phi\|_T<\infty$. The integral of
$\Phi$ is denoted by
$$\Phi \cdot \cw = \int_0^T \Phi_t\, d\cw_t$$
and the isometry property (\ref{iso}) is preserved for such
processes:
$$E(\|\Phi\punt \cw \|_H^2) = \| \Phi\|_T^2.$$
The details of this construction can be found in \cite[Chapter 4]{dz}.

   Let us conclude this section by providing a representation of the stochastic integral $\Phi\cdot \cw$ in terms of ordinary It\^o integrals of real-valued processes. Indeed, observe first that the expansion (\ref{wienerrepr}) can be rewritten in the form
\beq\label{3.10a}
   \cw_t=\sum_{j=1}^\infty \beta_j(t) \tilde e_j,
\eeq
where $(\tilde e_j)_j$ is defined in \eqref{e3.8}.

\begin{proposition}\label{prop3}
Let $(f_k)_k$ be a complete orthonormal system in the Hilbert space $H$. Assume that
$\Phi=\{\Phi_t,\ t\in [0,T]\}$ is any $L_2^0$-valued predictable process such that
$\|\Phi\|_T<\infty$. Then
\beq
\int_0^T \Phi_t\,
d\cw_t=\sum_{k=1}^\infty \left( \sum_{j=1}^\infty \int_0^T \langle \Phi_t(\tilde
e_j),f_k\rangle_H\, d\beta_j(t)\right) f_k.
\label{3}
\eeq
\end{proposition}

\noindent \emph{Proof.} First of all, we will prove that, under the standing hypotheses,
the right-hand side of (\ref{3}) is a well-defined element in $L^2(\Om;H)$. For this, we
will check that
$$
   E \left[\sum_{k=1}^\infty \left( \sum_{j=1}^\infty \int_0^T \langle \Phi_t(\tilde
e_j),f_k\rangle_H\, d\beta_j(t)\right)^2 \right]= \|\Phi\|_T^2,
$$
where the right-hand side is finite, by assumption.

Since $(\beta_j)_j$ is a family of independent standard Brownian motions,
$$
E \left[\sum_{k=1}^\infty \left( \sum_{j=1}^\infty \int_0^T \langle \Phi_t(\tilde
e_j),f_k\rangle_H\, d\beta_j(t)\right)^2 \right]
 = \sum_{k,j=1}^\infty
E\left[\left( \int_0^T \langle \Phi_t(\tilde e_j),f_k\rangle_H\, d\beta_j(t)\right)^2\right],
$$
and the right-hand side is equal to
$$
   \sum_{k,j=1}^\infty  \int_0^T E[\langle \Phi_t(\tilde e_j),f_k\rangle_H^2] \, dt
 =E \left[\int_0^T \sum_{j=1}^\infty \|\Phi_t(\tilde e_j)\|^2_H \, dt\right]
 = E \left[\int_0^T \|\Phi_t\|^2_{L^0_2} \, dt \right],
$$
and the last term is equal to $\|\Phi\|_T^2$. Hence, the series on the right-hand side of
(\ref{3}) defines an element in $L^2(\Om;H)$ and its norm is given by $\|\Phi\|_T$.
Therefore, by the isometry property of the stochastic integral (see (\ref{iso})), in order
to prove equality (\ref{3}), we only need to check this equality for simple processes. Namely, assume that $\Phi$ is of the form $\Phi_t=\Phi_0
1_{(a,b]}(t)$, where $\Phi_0$ is a $\tf_a$-measurable $L(V,H)$-valued random variable and
$0\leq a <b\leq T$. Then, by (\ref{simple_int}),
$$\int_0^T \Phi_t\, d\cw_t = \Phi_0 (\cw_b-\cw_a).$$
On the other hand,
$$
\sum_{k=1}^\infty \left( \sum_{j=1}^\infty \int_0^T \langle \Phi_t(\tilde
e_j),f_k\rangle_H\, d\beta_j(t)\right) f_k  = \sum_{k,j=1}^\infty
\langle \Phi_0(\tilde e_j),f_k\rangle_H (\beta_j(b)-\beta_j(a)) f_k,
$$
and the right-hand side is equal to
$$
   \sum_{j=1}^\infty  (\beta_j(b)-\beta_j(a))\, \Phi_0(\tilde e_j)
  = \Phi_0 \left( \sum_{j=1}^\infty   (\beta_j(b)-\beta_j(a)) \tilde e_j \right)
  = \Phi_0 (\cw_b-\cw_a),
$$
where the last equality follows from (\ref{3.10a}).
The proof is complete. \qed


\subsection{The case where $H=\R$}
\label{cylinf}

We consider a cylindrical Wiener process $W$ on some separable
Hilbert space $V$ with covariance $Q$, such that $\mbox{Tr } Q
<+\infty$. By Theorem \ref{t1}, $W$ is associated to a $V$-valued
$Q$-Wiener process $\cw$. We shall check that the
stochastic integral with respect to $W$, constructed in Section
\ref{integral}, is equal to an integral with respect to $\cw$,
constructed in \cite{dz} and sketched in Section \ref{daprato}, when the Hilbert space $H$ in which the integral takes its values is $H=\re$.

In Section \ref{integral}, we defined the Hilbert space $V_Q$ and the stochastic integral
$$g\cdot W = \int_0^T g_s\, dW_s,$$
for any predictable stochastic process $g\in L^2(\Om\times [0,T];
V_Q)$, with the isometry property
$$
E\left((g\cdot W)^2\right) = E\left( \int_0^T \|g_s\|^2_{V_Q}\, ds\right)^2.
$$
For any $s\in [0,T]$ and $g\in L^2(\Omega\times [0,T]; V_Q)$, we
define an operator $\Phi_s^g :V\rightarrow \re$ by \beq
   \Phi_s^g (\eta):= \langle g_s, \eta\rangle_{V},\qquad \eta\in V.\label{11}
\eeq
We denote by $L_2^0$ the set $L_2(V_0,H)$, with $V_0= Q^{1/2}(V)$ and $H=\re$.

\begin{proposition}\label{equiv}
Under the above assumptions, $\Phi^g =\{\Phi_s^g, \;s\in[0,T]\}$ defines a predictable
process with values in $L_2^0=L_2(V_0,\re)$, such that \beq E\left( \int_0^T
\|\Phi_s^g\|^2_{L_2^0}\, ds\right) = E\left( \int_0^T \|g_s\|^2_{V_Q}\,
ds\right).\label{equiv1} \eeq Therefore, the stochastic integral of $\Phi^g$ with respect
to $\cw$ can be defined as in Section \ref{daprato} and in fact,
\beq \int_0^T \Phi_s^g\, d\cw_s = \int_0^T g_s\, dW_s.\label{equiv2}
\eeq
\end{proposition}

\noindent\emph{Proof.} We first check (\ref{equiv1}).
Let $e_j$ be as in \eqref{wienerrepr}, $\tilde e_j$ be as in \eqref{e3.8} and $v_j$ be as in \eqref{e3.6}, so that $\tilde e_j = Q(v_j)$.
By \eqref{L02norm} with $H=\re$, and by \eqref{11},
$$
\|\Phi_s^g\|^2_2  = \sum_{j=1}^\infty  \langle g_s,\tilde e_j\rangle_V ^2
= \sum_{j=1}^\infty  \langle g_s,Q v_j \rangle_V^2
 = \sum_{j=1}^\infty \langle g_s, v_j \rangle_{V_Q}^2
= \|g_s\|^2_{V_Q}.
$$
We conclude that (\ref{equiv1}) holds. We note for later reference that this equality $\|\Phi_s^g\|_2 = \|g_s\|_{V_Q}$ remains valid even if Tr~$Q = + \infty$.

   Since, by hypothesis, the right hand-side of
(\ref{equiv1}) is finite, we deduce that $\Phi^g$ is a square
integrable process with values in $L_2^0$ and the
stochastic integral $\int_0^T \Phi^g_s\, d\cw_s$ is well-defined.

It remains to prove (\ref{equiv2}). For this, we apply Proposition \ref{prop3} in the
following situation: $H=\re$, with one basis vector $f_k=1$, $\Phi$ is defined in (\ref{11}), and
the sequence of independent standard Brownian motions in \eqref{3} is given by $\beta_j(t)=W_t(v_j)$.
Therefore,
$$
\int_0^T \Phi^g_t\, d\cw_t  = \sum_{j=1}^\infty \int_0^T \Phi^g_t(\tilde e_j) \,
d\beta_j(t),
$$
and the right-hand side is equal to
$$
   \sum_{j=1}^\infty \int_0^T \langle g_t,\tilde e_j\rangle_{V} \, dW_t(v_j)
   = \sum_{j=1}^\infty \int_0^T \langle g_t,v_j\rangle_{V_Q} \, dW_t(v_j)
  = \int_0^T g_t \,dW_t.
$$
This completes the proof.  \qed

\subsection{The case $\mbox{Tr } Q=+\infty$}
\label{cylwin}

In Proposition \ref{spcyl}, we showed that the covariance operator of the standard cylindrical Wiener process
$\{W_t(g),\; t\geq 0,\; g\in U\}$
associated with the spatially homogeneous noise that we considered in Section \ref{prel} is $Q = \mbox{Id}_U$, which
implies that $\mbox{Tr } Q=+\infty$. Therefore, we cannot make use of Proposition \ref{equiv} since, in this case,
there is no \emph{$Q$-Wiener process} associated to $W$. However,
there is the related notion of \emph{cylindrical $Q$-Wiener process}, which we now
define.

   Let $(V,\,\Vert\cdot\Vert_V)$ be a Hilbert space. Let $Q$ be a symmetric non-negative definite and bounded operator on $V$, possibly such that $\mbox{Tr } Q=+\infty$. Let $(e_j)_j$ be an orthonormal basis of $V$ that consists of eigenvectors of $Q$ with corresponding eigenvalues $\lambda_j$, $j \in \IN^*$. Define $V_0 = Q^{1/2} (V)$ as in Section \ref{daprato}.

   It is always possible to find a Hilbert space $V_1$ and a bounded linear injective
operator $J:(V,\|\cdot\|_V)\rightarrow (V_1,\|\cdot\|_{V_1})$ such that the restriction
$J_0=J_{|_{V_0}}:(V_0,\|\cdot\|_{V_0})\rightarrow (V_1,\|\cdot\|_{V_1})$ is
Hilbert-Schmidt. Indeed, as explained in \cite[Remark 2.5.1]{pr}, we may choose $V_1=V$,
$\langle \cdot,\cdot \rangle_{V_1} = \langle \cdot,\cdot \rangle_V$, $\al_k \in
(0,\infty)$ for all $k\geq 1$ such that $\sum_{k=1}^\infty \al_k^2 < +\infty$, and define
$J:V\rightarrow V$ by \beq J(h):= \sum_{k=1}^\infty \al_ k \langle h,e_k\rangle_V \;
e_k,\qquad h\in V, \label{4} \eeq where $(e_k)_k$ is an orthonormal basis of $V$. Then,
for $g\in V_0$, $g=\sum_{k=1}^\infty \langle g,\tilde e_k\rangle_{V_0} \; \tilde e_k$,
where $\tilde e_k= Q^{1/2} (e_k)$, $k\geq 1$, we have
$$J_0(g)= \sum_{k=1}^\infty \al_ k \langle g,\tilde e_k\rangle_{V_0} \sqrt{\lam_k}\; e_k=
\sum_{k=1}^\infty \al_ k \langle g,\tilde e_k\rangle_{V_0} \; \tilde e_k,$$ and so
$J_0:(V_0,\|\cdot\|_{V_0})\rightarrow (V,\|\cdot\|_V)$ is clearly Hilbert-Schmidt.

   As an operator between Hilbert spaces, from $V_0$ to $V_1$, $J_0$ has an adjoint $J_0^*: V_1 \to V_0$. However, if we consider $V_0$ and $V_1$ as Banach spaces, it is more common to consider the adjoint $\tilde J_0^*: V_1^* \to V_0^*$.

\begin{proposition}[{\cite[Proposition 4.11]{dz} and \cite[Proposition
2.5.2]{pr}}]\hspace{0.2cm}

    1.~Define $Q_1 = J_0 J_0^*:V_1= \mbox{\rm Im }J_0 \rightarrow V_1$. $Q_1$ is symmetric
    (self-adjoint), non-negative definite and $\mbox{\rm Tr }Q_1<+\infty$.

  2.~Let $\tilde e_j=Q^{1/2}(e_j)$, where $(e_j)_j$ is a complete orthonormal basis in
  $V$,  and let $(\beta_j)_j$ be a
family of independent  real-valued standard Brownian motions. Then
\beq \cw_t :=
\sum_{j=1}^\infty \beta_j(t) J_0(\tilde e_j),\qquad t\geq 0,
\label{10}
\eeq
is a $Q_1$-Wiener process in $V_1$.

   3.~Let $I:V_0\rightarrow V_0^*$ be the one-to-one mapping which identifies $V_0$ with its dual $V_0^*$, and consider the following diagram:
$$V_1^* \stackrel{\tilde J_0^*}{\longrightarrow} V_0^* 
\stackrel{I^{-1}}{\longrightarrow} V_0 \stackrel{J_0}{\longrightarrow} V_1.$$
Then, for all $s,t\geq 0$ and $h_1,h_2\in V_1^*$,
\beq \label{tif}
 E\left( \langle h_1,\cw_s\rangle_1\, \langle h_2,\cw_t\rangle_1 \right)=
(t\land s) \left\langle (I^{-1}\circ \tilde J^*_0) (h_1),(I^{-1}\circ \tilde J^*_0) (h_2)\right\rangle_{V_0},
\eeq
where $\langle \cdot,\cdot\rangle_1$ denotes the dual form on $V_1^*\times V_1$.

   4.~$\mbox{Im } Q_1^{1/2} = \mbox{Im } J_0$ and
$$\|h\|_0 = \|Q_1^{-1/2}  J_0 (h)\|_{V_1} = \| J_0(h)\|_{Q_1^{1/2}(V_1)}, \qquad  h\in V_0,$$
where $Q_1^{-1/2}$ denotes the pseudo-inverse of $Q_1^{1/2}$. Thus, $J_0: V_0 \rightarrow
Q_1^{1/2}(V_1)$ is an isometry.
\label{p1}
\end{proposition}

\begin{remark}
   (a) Part 3 in the Proposition's statement is commonly abbreviated in the following formal form (see, for instance, \cite[Proposition 1.1]{pz1}):
for all $s,t\geq 0$ and $h_1,h_2\in V_1^*$,
$$E\left( \langle h_1,\cw_s\rangle_1\, \langle h_2,\cw_t\rangle_1 \right)=
(t\land s) \left\langle h_1,h_2\right\rangle_{V_0}.$$

   (b) The $Q_1$-Wiener process $\{\cw_t,\; t\geq 0\}$ obtained in Proposition \ref{p1} is usually also called a
{\em cylindrical $Q$-Wiener process.} As it is pointed out in \cite[p.98]{dz}, if  $\mbox{\rm Tr
} Q<+\infty$, then we can take $\al_k=1$ in (\ref{4}), so $V_1=V$ and $J = \mbox{\rm Id}_V$,
and we get the classical concept of $Q$-Wiener process. In this case, one can take $V_0^* = V_Q$, $I^{-1} = Q\vert_{V_Q}$ and the equality \eqref{tif} reduces to
$$
   E\left( \langle h_1,\cw_s\rangle_1\, \langle h_2,\cw_t\rangle_1 \right)=
(t\land s) \left\langle Q h_1, h_2 \right\rangle_V.
$$
\label{remcyl}
\end{remark}

\noindent\emph{Proof of Proposition \ref{p1}.} Statement 1.~follows from \eqref{adjoint} and the fact that $J_0$ is Hilbert-Schmidt.
Concerning 2., we observe that for $h\in V_1$,
$$
E(\langle \cw_t,h\rangle_{V_1}^2)  = E\left( \left( \sum_{j=1}^\infty \beta_j(t) \langle
J_0(\tilde e_j),h\rangle_{V_1} \right)^2 \right),
$$
and the right-hand side is equal to
\begin{align*}
t \sum_{j=1}^\infty \langle J_0(\tilde e_j),h\rangle^2_{V_1}
 &= t \sum_{j=1}^\infty \langle \tilde e_j, J_0^*(h)\rangle^2_{V_0}
= t\, \|J_0^*(h)\|^2_{V_0}\\
&= t\, \langle J_0^*(h), J_0^*(h)\rangle_{V_0}
=t \, \langle J_0 J_0^*(h), h\rangle_{V_1}.
\end{align*}

Let us prove now part 3. For the sake of clarity, we will prove the statement for $s=t$ and $h_1=h_2$. Hence, let $t\geq 0$ and $h \in V_1^*$.
We denote by $\langle \cdot,\cdot \rangle_0$ the dual form on $V_0^*\times V_0$. Then, by (\ref{10}), the relation between $J_0$ and $J_0^*$, and the properties of $I$ and the family $(\beta_j)_j$, we obtain
$$
E(\langle h,\cw_t \rangle_1^2) = E\left( \left\langle h, \sum_{j=1}^\infty \beta_j(t) J_0(\tilde e_j)\right\rangle_1^2 \right),
$$
and the right-hand side is equal to
\begin{align*}
    t \sum_{j=1}^\infty \langle h,J_0(\tilde e_j)\rangle_1^2
& = t \sum_{j=1}^\infty \langle \tilde J_0^*(h),\tilde e_j\rangle_0^2
 = t \sum_{j=1}^\infty \langle (I^{-1}\circ \tilde J_0^*)(h),\tilde e_j\rangle_{V_0}^2 \\
&= t\, \| (I^{-1}\circ \tilde J_0^*)(h) \|^2_{V_0}.
\end{align*}
For 4., we refer the reader to  \cite[Proposition 2.5.2]{pr}. \qed

\vspace{0.25cm}

Let $\{\cw_t,\; t\geq 0\}$ be as in \eqref{10}. A predictable stochastic process
$\{\Phi_t,\; t\in [0,T]\}$ will be integrable with respect to $\cw$ if it takes values in
$L_2(Q_1^{1/2}(V_1),H)$ and
$$E\left( \int_0^T \|\Phi_t\|^2_{L_2(Q_1^{1/2}(V_1),H)}\, dt \right) <
+\infty.$$
By part 4 of Proposition \ref{p1}, we have
$$
   \Phi \in L_2^0 = L_2(V_0,H)\quad \Longleftrightarrow\quad \Phi \circ J_0^{-1} \in  L_2 (Q_1^{1/2}(V_1),H).
$$

\begin{definition}
   For any square integrable
predictable process $\Phi$ with values in $L_2^0$ such that
$$E\left( \int_0^T \|\Phi_t\|^2_{L_2^0}\, dt \right) <+\infty,$$
the $H$-valued stochastic integral $\Phi \cdot \cw$ is defined by
$$\int_0^T \Phi_s\, d\cw_s := \int_0^T \Phi_s \circ J_0^{-1}\, d\cw_s.$$
\end{definition}

   We note that the class of integrable processes with respect to $\cw$  does not depend on
the choice of $V_1$.

 We now relate this notion of stochastic integral with the stochastic integral with respect to the cylindrical Wiener process of Section \ref{integral}. Let $\{W_t, \; t\in [0,T]\}$ be a cylindrical Wiener process with covariance $Q$ on the
Hilbert space $V$, and let $g\in L^2(\Om\times [0,T];V_Q)$ be a predictable process, so that $g\cdot W$ is well defined as
in Section \ref{integral}. By Proposition \ref{p1}, we can consider the cylindrical $Q$-Wiener process $\{\cw_t,\;
t\in [0,T]\}$ defined by
\beq\label{QWP}
   \cw_t = \sum_{j=1}^\infty \beta_j(t) J_0(\tilde e_j)
\eeq
as in formula (\ref{10}) with $\beta_j(t)=W_t(v_j)$, where
$v_j=Q^{-1/2}(e_j)$, $\tilde e_j  = Q^{1/2}(e_j)$ and $(e_j)_j$ denotes a complete orthonormal basis in $V$ consisting of eigenvalues of $Q$, so that
$(v_j)_j$ is a complete orthonormal basis in $V_Q$. This process takes values in some Hilbert space $V_1$.

For $g\in L^2(\Om\times [0,T];V_Q)$, we define, as in (\ref{11}), the operator
$$\Phi^g_s(\eta)=\langle g_s,\eta\rangle_V, \qquad \eta\in V,$$
which takes values in $H=\re$. Recall that $V_0 = Q^{1/2}(V)$ and $V_Q=Q^{-1/2}(V)$.

\begin{proposition}\label{prop24}
The process $\{\Phi_s^g,\, s\in [0,T]\}$ defines a predictable process with values in
$L_2(V_0,\re)$, such that
$$E\left( \int_0^T \|\Phi_s^g\|^2_2\, ds\right)=E\left( \int_0^T \|g_s\|^2_{V_Q}\, ds \right),$$ and
$$ \int_0^T \Phi_s^g\, d\cw_s= \int_0^T g_s\, dW_s.$$
\end{proposition}

\noindent{\emph{Proof.}} First, we will prove that $\Phi_s^g\in L_2(V_0,\re)$, for $s\in [0,T]$. As in the first part of the proof of Proposition \ref{equiv}, $\|\Phi_s^g\|_2 = \|g_s\|_{V_Q}$. This gives the equality of expectations in the statement of the proposition, and the right-hand side is finite by assumption.

Concerning the equality of integrals, we note that by definition,
$$\int_0^T \Phi_s\, d\cw_s := \int_0^T \Phi_s \circ J_0^{-1}\, d\cw_s,$$
where the right-hand side is defined using the finite-trace approach of Section
\ref{daprato}. We note that by Proposition \ref{p1}, part 4, $(J_0(\tilde e_j))_j$ is a
complete orthonormal basis of $Q_1^{1/2}(V_1)$.

According to Proposition \ref{prop3} with $H=\R$, a single basis element $f_k=1$ of $H$,
$\beta_j(s)=W_s(v_j)$, and $\tilde e_j$ there replaced by $J_0(\tilde e_j)$, formula
(\ref{3}) becomes
$$\int_0^T \Phi^g_s \circ J_0^{-1}\, d\cw_s = \sum_{j=1}^\infty \int_0^T \Phi_s^g \circ
J_0^{-1} (J_0(\tilde e_j))\, d\beta_j(s),$$ and the right-hand side is equal to
\begin{align*}
\sum_{j=1}^\infty \int_0^T \Phi_s^g (\tilde e_j)\, dW_s(v_j) &= \sum_{j=1}^\infty
\int_0^T \langle g_s,\tilde e_j\rangle_V \, dW_s(v_j)\\
&=\sum_{j=1}^\infty \int_0^T \langle g_s,v_j\rangle_{V_Q} \, dW_s(v_j).
\end{align*}
The last expression is equal to $g\cdot W$. \qed

\begin{remark}\label{r1}
Proposition \ref{prop24} allows us in particular to associate the spatially homogeneous noise of Section \ref{prel}, viewed as a cylindrical Wiener process with covariance {\rm Id}$_U$ in Proposition \ref{spcyl}, with a cylindrical $Q$-Wiener process as in Proposition \ref{p1}, with $Q=\mbox{\rm Id}_U$, on the Hilbert space $U$ of Section \ref{prel}, and to relate the associated stochastic integrals.
\end{remark}


\section{Spde's driven by a spatially homogeneous noise}
\label{spde}

This section is devoted to presenting a class of spde's in $\red$ driven by a spatially homogeneous noise. In Section \ref{rva}, we present the real-valued approach using the notion of a mild random field solution of
the equation. Section \ref{examples} gives two examples: the
stochastic heat equation in any spatial dimension and the stochastic
wave dimension in spatial dimensions $d=1,2,3$. In Section \ref{eu}, we establish an existence and uniqueness result which extends a theorem of \cite{dalang}. In Section \ref{homo-inf}, we present the infinite-dimensional formulation of these spde's. In Section \ref{relation}, we examine the relationship between these two formulations, and conclude that they are equivalent (see Proposition \ref{rdequiv}). In Section \ref{relation-dm}, we examine the relationship with the approach of \cite{dm}.

We are interested in the following class of non-linear spde's:
\begin{equation}
L u(t,x) = \sigma(u(t,x)) \dot{W} (t,x) + b(u(t,x)),
\label{2.1}
\end{equation}
$t \geq 0$, $x\in \mathbb{R}^{d}$, where $L$ denotes a general
second order partial differential operator with constant coefficients, with appropriate initial
conditions. The coefficients $\sigma$ and $b$ are real-valued
functions and $\dot{W} (t,x)$ is the formal notation for the
Gaussian random perturbation described at the beginning of Section
\ref{prel}.

If $L$ is first order in time, such as the heat operator
$L=\frac{\partial}{\partial t}-\Delta$, where $\Delta$ denotes the Laplacian operator on $\mathbb{R}^d$, then we impose initial
conditions of the form \beq u(0,x)=u_0(x) \qquad x\in \red,
\label{2.1.1} \eeq for some Borel function $u_0:\red \rightarrow \re$. If $L$ is second order in time, such as the wave operator
$L=\frac{\partial^2}{\partial t^2}-\Delta$, then we have to impose two
initial conditions: \beq u(0,x)=u_0(x), \quad \frac{\partial
u}{\partial t}(0,x)=v_0(x),\qquad x\in \red, \label{2.1.2} \eeq for
some Borel functions $u_0,v_0:\red \rightarrow \re$.


\subsection{The random field approach}
\label{rva}

We now describe the notion of \emph{mild random field} solution to equation (\ref{2.1}). Recall that we are given a filtered probability space $(\Om,\tf,(\tf_t),P)$, where $(\tf_t)_t$ is the filtration generated by the standard cylindrical Wiener process $W$ of Proposition \ref{spcyl}, and we fix a time
horizon $T>0$. A real-valued adapted stochastic process $\{u(t,x),\, (t,x)\in [0,T]\times \red\}$ is a \textit{mild random field} solution of (\ref{2.1}) if the following stochastic integral equation is satisfied:
\begin{align}
u(t,x)=& I_0(t,x) + \int_0^t \int_{\mathbb{R}^{d}} \Gamma(t-s,x-y) \sigma(u(s,y))\, W(ds,dy) \nonumber\\
&+ \int_0^t ds \int_{\mathbb{R}^{d}} \Gamma(s,dy)\, b(u(t-s,x-y)), \qquad a.s.,
\label{2.2}
\end{align}
for all $(t,x) \in [0,T] \times \red$.
In (\ref{2.2}), $\Gamma$ denotes the fundamental solution
associated to $L$ and $I_0(t,x)$ is the contribution of the initial
conditions, which we define below. The stochastic integral on the
right hand-side of (\ref{2.2}) is as defined in Section \ref{homo}. In particular, we need to assume that for any $(t,x)$,
the fundamental solution $\Gam(t-\punt,x-\star )$ satisfies Hypothesis \ref{hypA},
and  to require that $s\mapsto \Gam(t-s,x-\star ) \sig(u(s,\star ))$, $s\in [0,t]$, defines a
predictable process taking values in the space $U$ of Section \ref{prel} such that
$$E\left( \int_0^t \|\Gam(t-s,x-\star ) \sig(u(s,\star ))\|^2_U\, ds \right)
<+\infty$$
(see Sections \ref{prel} and \ref{distri}). As we will make explicit in Section \ref{eu}, these assumptions will be satisfied under certain regularity
assumptions on the coefficients $b$ and $\sig$ (see Theorem
\ref{existence}).

The last integral on the
right-hand side of (\ref{2.2}) is considered in the pathwise sense,
and we use the notation ``$\Gam(s,dy)$'' because we
will assume that $\Gam(s)$ is a measure on $\red$. Concerning the term $I_0(t,x)$, if $L$ is a parabolic-type operator
and we consider the initial condition (\ref{2.1.1}), then
\beq I_0(t,x)=
\left( \Gam(t)
* u_0\right)(x)=\int_{\red} u_0(x-y)\, \Gam(t,dy). \label{2.3} \eeq On the other hand, in
the case where $L$ is second order in time with initial values
(\ref{2.1.2}),
\begin{align}
I_0(t,x) & = \left( \Gam(t) * v_0\right)(x) + \frac{\partial}{\partial t} \left( \Gam(t) * u_0\right)(x)\nonumber \\
& =\int_{\red} v_0(x-y)\, \Gam(t,dy) + \frac{\partial}{\partial t} \left(\int_{\red} u_0(x-y)\, \Gam(t,dy)\right).
\label{2.4}
\end{align}


\subsection{Examples: stochastic heat and wave equations}
\label{examples}

   In the case of the stochastic heat equation in any space dimension $d\geq 1$ and the stochastic wave equation in
dimensions $d=1,2,3$, following \cite[Section 3]{dalang} (see also \cite[Examples 4.2 and 4.3]{nualartquer}), the fundamental solutions are well-known and the conditions in Hypothesis \ref{hypA} can be made explicit.

  Indeed, let $\Gam$ be the fundamental solution of the heat
equation in $\red$, $d\geq 1$, so that
$$\Gam(t,x)=(4\pi t)^{-d/2} \exp\left(-\frac{|x|^2}{4t}\right).$$
In particular, we have $\mathcal{F} \Gamma (t)(\xi)=\exp(-4\pi^2
t|\xi|^2)$, $\xi\in \red$, and, because
$$\int_0^T \exp(-4\pi^2 t|\xi|^2)\, dt =\frac{1}{4\pi^2 |\xi|^2} (1-\exp(-4\pi^2 T|\xi|^2)),$$
we conclude that condition (\ref{mu-L}) in  Hypothesis \ref{hypA}
holds if and only if
\beq\int_{\red}
\frac{\mu(d\xi)}{1+|\xi|^2}<+\infty.
\label{18}
\eeq


   Now let $\Gam_d$ be the fundamental solution of the wave equation in $\red$, with $d=1,2,3$. This restriction on the space
dimension is due to the fact that the fundamental solution in $\red$ with $d>3$ is no longer a non-negative distribution (for results on the stochastic wave equation in spatial dimension $d>3$, we refer the reader to \cite{CD}: see Remark \ref{conus-dalang}). It is well known (see \cite[Chapter 5]{folland}) that
$$
\Gamma_1(t,x)=\frac{1}{2} 1_{\{|x|<t\}}, \quad
\Gamma_2 (t,x) = \frac{1}{2\pi} (t^2-|x|^2)_+^{-1/2}, \quad
\Gamma_3(t)(dx) = \frac{1}{4\pi t}\sigma_t(dx),
$$
where $\sig_t$ denotes the uniform surface measure on the
three-dimensional sphere of radius $t$, with total mass $4 \pi t^2$. This implies that, for each
$t$, $\Gamma_d(t)$ has compact support. Furthermore, for all
dimensions $d\geq 1$, the Fourier
transform of $\Gamma_d(t)$ is
$$\mathcal{F} \Gamma_d (t)(\xi)= \frac{\sin (2\pi t |\xi|)}{2\pi |\xi|}.$$
Elementary estimates show that there are positive constants $c_1$
and $c_2$ depending on $T>0$ such that
$$\frac{c_1}{1+|\xi|^2}\leq \int_0^T \frac{ \sin^2 (2\pi t\xi)}{4\pi^2 |\xi|^2} dt \leq \frac{c_2}{1+|\xi|^2}.$$
Therefore, $\Gamma_d$ satisfies condition (\ref{mu-L}) if and only
if (\ref{18}) holds.

For $d=1$, $I_0(t,x)$ is given by the so-called d'Alembert's formula
(see, for instance, \cite[p.68]{evans}): \beq I^1_0(t,x)=
\frac{1}{2} \left[ u_0(x+t)+u_0(x-t)\right] +
\frac{1}{2}\int_{x-t}^{x+t} v_0(y)\, dy, \qquad x\in \re. \label{2.4.2}
\eeq For $d=2$ (see \cite[p.74]{evans}), 
$$I^2_0(t,x)= \frac{1}{2\pi t} \int_{|x-y|<t} \frac{u_0(y + t v_0 )+ \nabla u_0(y)\cdot (x-y)}{(t^2- |x-y|^2)^{1/2}}\, dy, \qquad x\in \re^2.$$
Finally, for $d=3$ (see
\cite[p.77]{evans}), for $x\in \re^3$,
  \beq I^3_0(t,x)= \frac{1}{4\pi
t^2} \int_{\re^3} \left( tv_0(x-y) + u_0(x-y) +  \nabla
u_0(x-y)\cdot y \right)\, \sig_t(dy). \label{2.4.3} \eeq
It is important to remark that in the above formulas, we have
implicitly  assumed that all integrals that appear are well
defined. Indeed, in Lemma \ref{lemma2}
below, we will exhibit sufficient conditions on $u_0$ and $v_0$
under which such integrals exist and are uniformly bounded with respect to $t$ and $x$.

\subsection{Random field solutions with arbitrary initial conditions}
\label{eu}



The aim of this section is to prove the existence and uniqueness of a mild random field solution to the stochastic integral equation \eqref{2.2}.

   We are interested in solutions that are $L^p$-bounded, as in \eqref{6.25a} below, and $L^2$-continuous. This is only possible under certain assumptions on the initial conditions. In particular, the initial conditions will have to be such that the following hypothesis is satisfied.
   
\begin{hypothesis}\label{hyp0}
$(t,x) \mapsto I_0(t,x)$ is continuous and 
$\sup_{(t,x)\in[0,T]\times \red} |I_0(t,x)|<+\infty$.
\end{hypothesis}

For the particular case of the heat equation in any spatial
dimension and the wave equation with $d \in \{1,2,3\}$, sufficient conditions for Hypothesis \ref{hyp0} to hold are given in the next lemma.

\begin{lemma}\label{lemma2}
Consider the following two sets of hypotheses:
\begin{itemize}
    \item[(i)] {\em Heat equation.} $u_0:\red\rightarrow \re$ is measurable and bounded.
    \item[(ii)] {\em Wave equation.} When $d=1$, $u_0$ is bounded and continuous, and $v_0$ is bounded and measurable. When $d=2$, $u_0 \in C^1(\re^2)$ and there is $q_0 \in \, ]2,\infty]$ such that $u_0,\nabla u_0, v_0$ all belong to $L^{q_0}(\re^2)$. When $d=3$, $u_0 \in C^1(\re^3)$, $u_0$ and $\nabla u_0$ are bounded, and $v_0$ is bounded and continuous.
\end{itemize}
Then under condition (i) or (ii), Hypothesis \ref{hyp0} is
satisfied.
\end{lemma}

\noindent\textit{Proof.} Assume first that $L$ is the heat operator on
$\red$, $d\geq 1$,  with initial condition $u_0$ satisfying (i).
Then, by (\ref{2.3}), 
\begin{align*}
   \sup_{(t,x)\in [0,T]\times \red} |I_0(t,x)| &\leq \|u_0\|_{\infty} \sup_{t\in\, ]0,T]}  \int_{\red} (2\pi t)^{-d/2}\exp\left(-\frac{|y|^2}{2t}\right) dy \\
    &=\|u_0\|_{\infty}<+\infty.
\end{align*}

   Secondly, assume that  $L$ is the wave operator on $\red$,
$d=1,2,3$, and that condition (ii) is satisfied. We make
explicit the dependence on the space dimension by denoting $I_0^d(t,x)$, $d=1,2,3$, the term $I_0(t,x)$.

By (\ref{2.4.2}), if $d=1$ it is clear that
$$\sup_{(t,x)\in [0,T]\times \re}  |I_0^1(t,x)| \leq C (\|u_0\|_{\infty} + \|v_0\|_{\infty}).$$
To deal with the case $d=2$,  we refer to \cite[p.808--809]{MS}. In
this reference, the explicit formula $\Gam_2(t,x)=\frac{1}{2\pi}
(t^2-|x|^2)_+^{-1/2}$ was used to show that
$$\sup_{(t,x)\in [0,T]\times \re}  |I_0^2(t,x)| \leq C (\|u_0\|_{\infty} + \|\nabla u_0\|_{\infty} + \|v_0\|_{\infty}).$$
Finally, for the case $d=3$ we have, by (\ref{2.4.3}):
$$|I_0^3(t,x)| \leq C(\|v_0\|_{\infty} + \|u_0\|_{\infty}+ \|\nabla u_0\|_{\infty}) \sup_{s\in\, ]0,T]}
\frac{\sig_s(\re^3)}{s^2},$$ where $\sig_s$ denotes the uniform
surface measure on the three-dimensional sphere of radius $s$. In
particular, the total mass of $\sig_s$ is proportional to $s^2$ and, therefore, $I_0^3(t,x)$ is uniformly bounded with respect to $t\in [0,T]$ and $x\in \re^3$. 

   Finally, the continuity property of $(t,x) \mapsto I_0(t,x)$ follows from the hypotheses and the explicit formulas for $I_0(t,x)$ given in Section \ref{rva}. This concludes the proof. \qed
\vskip 16pt

   The next theorem discusses existence and uniqueness of mild random field solutions to equation \eqref{2.2}. Since this theorem covers rather general initial conditions, it is an extension of Theorem 13 in \cite{dalang}. Indeed, in this reference, only vanishing initial data could be considered, because of the spatially homogeneous covariance required for the process $Z$ in the construction of the stochastic integral used there for the wave equation when $d=3$ (see \cite[p.10 and Theorem 2]{dalang}). Of course, in the case of the stochastic wave equation in spatial dimensions $d=1,2$, there are
many results on existence and uniqueness of mild random field solutions with non-vanishing initial conditions: see for instance \cite{carmona,DF,MS,mueller}.

\begin{theorem}\label{existence}
Assume that Hypotheses \ref{hypA} and \ref{hyp0} are satisfied and that $\sig$ and $b$ are Lipschitz functions. Then there exists a unique mild random field solution $\{u(t,x),\, (t,x)\in [0,T]\times \red\}$ of equation (\ref{2.2}). Moreover, the process $u$ is $L^2$-continuous and for all $p\geq 1$,
\begin{equation}\label{6.25a}
\sup_{(t,x)\in [0,T]\times \red} E(|u(t,x)|^p) <+\infty.
\end{equation}
\end{theorem}

\noindent\textit{Proof.} The proof is similar to those of \cite[Theorem
1.2]{MS} and \cite[Theorem 13]{dalang}. We define the Picard
iteration scheme
\begin{align}
u^0(t,x)&=  I_0(t,x),\nonumber \\
u^{n+1}(t,x) &=  u^0(t,x) + \int_0^t \int_{\red} \Gam(t-s,x-y) \sig(u^n(s,y))\, W(ds,dy) \nonumber \\
&\qquad + \int_0^t \int_{\red} b(u^n(t-s,x-y))\,  \Gam(s,dy)\, ds,
\label{2.7}
\end{align}
for $n\geq 0$. We prove by induction on $n$ that the process
$\{u^n(t,x), \;(t,x)\in [0,T]\times \red\}$ is well defined and, for
 $p\geq 1$,
\beq \sup_{(t,x)\in [0,T]\times \red} E(|u^n(t,x)|^p) <+\infty,
\label{6.25} \eeq for every $n\geq 0$.

Notice that by Hypothesis \ref{hyp0}, the process $u^0$ is
locally bounded, and the Lipschitz property on $\sig$ yields
$$\sup_{(t,x)\in [0,T]\times \red} |\sig(u^0(t,x))|^p <+\infty.$$
By Proposition \ref{prop1}, this implies that the stochastic
integral 
$$
   \ci^0(t,x)=\int_0^t \int_{\red}
\Gam(t-s,x-y)\sig(u^0(s,y))\, W(ds,dy)
$$ 
is well-defined and
\begin{align}
E(|\ci^0(t,x)|^p)& \leq C \int_0^t ds \sup_{z\in \red}
\left(1+|u^0(s,z)|^p\right) \int_{\red} \mu(d\xi)\, |\tf \Gam
(t-s)(\xi)|^2\nonumber \\
& \leq C \sup_{(s,z)\in [0,T]\times \red} (1+|u^0(s,z)|^p) \int_0^T ds\, J(s), \label{6.3}
\end{align} 
where 
$$
   J(s)=\int_{\red} \mu(d\xi)\, |\tf \Gam(s)(\xi)|^2.
$$ 
In order to deal with the pathwise integral 
$$
   \cj^0(t,x)=\int_0^t ds \int_{\red} \Gam(s,dy)\,
b(u^0(t-s,x-y)),
$$
we apply H\"older's inequality with respect to the
finite measure $\Gam(s,dy)ds$ on $[0,T]\times \red$ and use the Lipschitz
property of $b$ : 
\beq 
   |\cj^0(t,x)|^p \leq C \int_0^t ds \int_{\red} \Gam(s,dy) \left( 1+|u^0(t-s,x-y)|^p\right). \label{6.4} 
\eeq 
The latter term is uniformly bounded with respect to $t$ and $x$.
Together with (\ref{6.3}), this implies that $\{u^1(t,x),\;(t,x)\in [0,T]\times \red\}$ is a well-defined measurable process. Further, by (\ref{6.3}), (\ref{6.4}) and Hypothesis \ref{hypA},
$$
   \sup_{(t,x)\in [0,T]\times \red} E(|u^1(t,x)|^p) <+\infty.
$$

   Consider now $n>1$ and assume that $\{u^n(t,x),\;(t,x)\in [0,T]\times
\red\}$ is a well-defined measurable process satisfying
(\ref{6.25}). Using the same arguments as above, one proves
that the integrals $\ci^{n+1}(t,x)$ and $\cj^{n+1}(t,x)$ exist, so that the process $u^{n+1}$ is well-defined and is uniformly bounded in $L^p(\Om)$. This proves \eqref{6.25}.

   The next step consists in showing that the bound
(\ref{6.25}) is uniform with respect to $n$, that is
\beq \sup_{n\geq 0}
\sup_{(t,x)\in [0,T]\times \red} E(|u^n(t,x)|^p)
<+\infty.\label{6.1} 
\eeq 
Indeed, the same kind of estimates as in
the first part of the proof show that for $n\geq 1$,
$$E(|u^{n+1}(t,x)|^p)\leq C\left(1+ \int_0^t ds \left( 1+\sup_{z\in
\red} E(|u^n(s,z)|^p)\right) (J(t-s)+1) \right),$$ 
We conclude that (\ref{6.1}) holds by the version of Gronwall's Lemma presented in \cite[Lemma 15]{dalang}.

Now we show that the sequence $(u^n(t,x))_{n\geq 1}$ converges in
$L^p(\Om)$. Following the same lines as in the proof of \cite[Theorem 13]{dalang}, let
$$M_n(t):= \sup_{(s,x)\in [0,t]\times \red}
E(|u^{n+1}(s,x)-u^n(s,x)|^p).$$ Using the Lipschitz property of $b$
and $\sig$, and applying the same arguments as above, we obtain the
estimate
$$M_n(t)\leq C \int_0^t ds\, M_{n-1}(s)(J(t-s)+1).$$
Hence, we apply again \cite[Lemma 15]{dalang} to
conclude that $(u^n(t,x))_{n\geq 1}$ converges uniformly in
$L^p(\Om)$ to a limit $u(t,x)$. The process $\{u(t,x),\ (t,x)\in
[0,T]\times \red\}$ has a measurable version that satisfies equation (\ref{2.2}). Indeed, let us sketch the calculations concerning the
stochastic integral term $\ci^n(t,x)$ of (\ref{2.7}): we will prove that
$$\lim_{n\rightarrow \infty} \sup_{(t,x)\in [0,T]\times \red}
E(|\ci^n(t,x)-\ci(t,x)|^p) =0.$$ By the Lipschitz property of
$\sig$, Proposition \ref{prop1} and Hypothesis \ref{hypA},
\begin{align*}
& E(|\ci^n(t,x)-\ci(t,x)|^p) \\
&\quad   \leq E\left(\Big\vert\int_0^T \int_{\red}
\Gam(t-s,x-y)[\sig(u^{n-1}(s,y))-\sig(u(s,y))]\, W(ds,dy)\Big\vert^p\right)\\
& \quad \leq C \int_0^T ds \sup_{z\in \red}
E(|u^{n-1}(s,z)-u(s,z)|^p)
\int_{\red} \mu(d\xi)\, |\tf \Gam(t-s)(\xi)|^2\\
& \quad \leq C \sup_{(s,z)\in [0,T]\times \red}
E(|u^{n-1}(s,z)-u(s,z)|^p)
\end{align*}
and this last term converges to zero, as $n$ tends to infinity. The
pathwise integral term can be studied in a similar manner.
Therefore, the process $u$ solves (\ref{2.2}). Finally, uniqueness of the solution can be checked by standard arguments. \qed

\subsection{Spatially homogeneous spde's in the infinite-dimensional setting}
\label{homo-inf}

Stochastic partial differential equations of the form (\ref{2.1}) on
$\red$ and driven by a spatially homogeneous Wiener process have
been studied, in the context of Da Prato and Zabczyk \cite{dz}, in a series of works: \cite{kz1,kz2,peszat,pz1,pz2}. The aim of this section is to sketch the formulation used in those papers,  focusing mostly on the one used by Peszat and Zabczyk in \cite{pz2}. Then, in Section \ref{relation}, we will compare their solution with the mild random field solution of Section \ref{rva}. 

   In \cite{pz2}, the stochastic wave equation with $d=1,2,3$ and the stochastic heat equation in any space dimension are considered. This meshes well with the case considered in Section \ref{rva}, in which the fundamental solution associated to the underlying differential operator is a non-negative distribution. However, we note that the stochastic wave equation in higher dimensions ($d>3$) can also be formulated and solved in the infinite-dimensional setting, but using a slightly different formulation (see \cite{peszat}).

To begin with, we notice that in \cite{pz2}---and, indeed, in the above mentioned companion papers---the slightly more general
spatially correlated noise described in Remark \ref{rem0} is used. More precisely, one considers a spatially homogeneous Wiener process $\{\cw^*_t,\; t\geq 0\}$ with values in the space $\cs'(\red)$ of tempered
distributions. If we denote by $\langle \punt, \punt\rangle$ the
usual duality action of $\cs'(\red)$ on $\cs(\red)$, this means that for all $\ffi\in \cs(\red)$, $\{\langle \cw^*_t,\ffi\rangle,\; t\in \re_+\}$ is a 
centered Gaussian process and there exists $\Lam\in \cs'(\red)$ such that for all $\ffi,\psi\in \cs(\red)$ and $s,t\in \re_+$,
$$E\left( \langle \cw^*_s, \varphi\rangle \langle \cw^*_t, \psi\rangle \right) = (s\wedge t)\langle \Lambda, \varphi * \tilde
\psi\rangle,$$
where $\tilde \psi(x)=\psi(-x)$. The Schwartz distribution $\Lam$ must be the Fourier transform of a symmetric and non-negative tempered measure $\mu$ on $\red$.



\begin{remark}
  In the particular case where $\Lam$ has a density $f$ satisfying the conditions of Section \ref{prel}, we
recover the covariance operator of the cylindrical Wiener process
$W$ on the Hilbert space $U$ defined in Proposition \ref{spcyl} (see
\eqref{1}):
\begin{align}
E\left( \langle \cw^*_s, \varphi\rangle \langle \cw^*_t, \psi\rangle \right) & = (s\wedge t) \int_{\red} dx  \int_{\red}
dy\, \vp(x) f(x-y) \psi(y) \nonumber \\
& = E\left( W_s(\ffi) W_t(\psi)\right).
\label{56}
\end{align}
\label{rem4.1}
\end{remark}

   For the sake of clarity in the exposition, we assume that the spatial correlation of the noise is given by $\Lam=f$,
as just described in Remark \ref{rem4.1}.


   Let $U$ be the Hilbert space defined in Section \ref{prel}, and let $U^*$ be the dual of $U$. The following characterization of $U^*$ is given in \cite[Proposition 1.2]{pz1}. Recall that $\tilde L^2
(\red, \mu)$ stands for the subspace of $L^2(\red, \mu)$ consisting of
all functions $\phi$ such that $\tilde \phi =\phi$.

\begin{lemma}\label{rdU'}
A distribution $g\in \cs'(\red)$ belongs to $U^*$ if and only if
there is $\phi\in \tilde L^2(\red, \mu)$ such that $g=\tf (\phi
\mu)$. Moreover, if $g_1 =\tf (\phi_1 \mu)$ and $g_2 =\tf (\phi_2
\mu)$, with $\phi_1,\phi_2 \in \tilde L^2(\red, \mu)$,  then
$$\langle g_1, g_2\rangle_{U^*} = \langle \phi_1, \phi_2\rangle_{\tilde L^2 (\red, \mu)}.$$
\end{lemma}

\begin{remark}\label{r2}
The previous lemma allows us to
determine the explicit form of the isometry $I:U\rightarrow U^*$. More precisely, as stated in Remark
\ref{rem1}, any element $g \in U$ can be written in the form
$g=\tf^{-1} \phi$, with $\phi\in \tilde L^2(\red; d\mu)$. Then, for
such $g$, $I(g)\in U^*$ is defined by
$$I(g)=\tf(\phi \mu).$$

\end{remark}

   Moreover, we have the following lemma whose proof is
straightforward. In this lemma, $\tilde \cs(\red)$ denotes the family of
functions $\vp\in \cs(\red)$ such that $\tilde\vp = \vp$.

\begin{lemma}
Let $\ffi\in U$ be such that $\ffi\in \tilde \cs(\red)$. Then $I(\ffi)=\ffi*f$.
\label{lemma3}
\end{lemma}

   As it has been explained in \cite[p.191]{pz1}
(see, in particular, Proposition 1.1 therein), $\cw^*$ may be regarded as a $U^*$-valued cylindrical $Q$-Wiener process with $Q = \mbox{ Id}_{U^*}$. More precisely, let $U_1^*$ be a
Hilbert space such that there exists a dense Hilbert-Schmidt embedding $J^*:U^*\rightarrow U^*_1$ (see Proposition \ref{p1}). Then
\beq\label{wstarrd}
   \cw^*_t = \sum_{j=1}^\infty \beta_j(t) J^*(e_j^*),
\eeq
where $(e_j^*)_j$ is a complete orthonormal basis in $U^*$, and the $\beta_j(t)$ are independent standard Brownian motions (note that $Q^{1/2} = Q^{-1/2}=\mbox{ Id}_{U^*}$). Therefore, we will be able to define Hilbert-space-valued stochastic integrals with respect to $\cw^*$, as has been described in Section \ref{cylwin}. Note that $U^*$ is sometimes called the reproducing kernel Hilbert space associated to $\cw^*$ (see, for instance, \cite[Section 2.2.2]{dz} or \cite[p.191]{pz1}).


In \cite{pz2}, mild solutions to the formal equation
(\ref{2.1}) are considered in a Hilbert space $H$ of the form $L^2_\vt = L^2(\red, \vt(x)dx)$, where
$\vt$ is a strictly positive even function such that $\vt(x)=e^{-|x|}$, for $|x|\geq 1$. Let us also denote by $H_\vt^1$ the weighted Sobolev space
which is the completion of $\cs(\red)$ with respect to the norm
$$\|\psi\|_{H_\vt^1} =\left( \int_{\red} \left[ |\psi(x)|^2 + |\nabla \psi(x)|^2 \right] \vt(x) dx \right)^{1/2}.$$

For simplicity, we will restrict ourself to equation (\ref{2.1}) when $L=\frac{\partial^2}{\partial t^2}- \Del$ with spatial dimensions $d=1,2,3$. Let $\Gam$ be the fundamental solution associated to $L$, let $u_0\in H^1_\vt$, $v_0\in L^2_\vt$, and fix a time horizon
$T>0$. By definition, a \textit{mild $L^2_\vt$-valued} solution of (\ref{2.1}) with
$L=\frac{\partial^2}{\partial t^2} -\Del$, is an $\tf_t$-adapted
process $\{u(t),\; t\in [0,T]\}$ with values in $L^2_\vt$ satisfying
\begin{align}
u(t)  = & \frac{\partial}{\partial t} \left( \Gam(t)*u_0 \right) +
\Gam(t)* v_0 + \int_0^t \Gam(t-s) * b(u(s))\, ds \nonumber \\
& + \int_0^t \Gam(t-s) * \sig(u(s))\, d\cw^*_s. \label{8}
\end{align}
The stochastic integral on the right-hand side of (\ref{8}) has to be defined. This requires interpreting the integrand $\Gam(t-s) * \sig(u(s))$ in the framework of Section \ref{cylwin}.


  Recall that, as in Section \ref{cylwin} and since $Q = \mbox{ Id}_{U^*}$ and so $U^* = (U^*)_0$, we will be able to
define the stochastic integral with respect to $\cw^*$ of any
predictable process $\Phi$ taking values in the space $L_2(U^*,H)$, where $H = L^2_\vt$. Therefore, it is necessary to interpret $\Gam(t-s) * \sig(u(s))$ as an element of $L_2(U^*,H)$.

Let $U^{*,0}$ be the dense subspace of $U^*$ consisting of all $g=\tf
(\phi \mu)$ with $\phi\in \tilde \cs(\red)$. According to
\cite[p.427]{pz2}, it holds that $U^{*,0}\subset \cac_b(\red)$, the
space of bounded and continuous functions on $\red$. For $u\in
L^2_\vt$ and $t>0$, define the following operator:
\beq
\ck(t,u)(\eta) = \Gam(t)*(u \eta),\qquad \eta\in U^{*,0}.
\label{994}
\eeq
Then it is shown in Lemma 3.3 of \cite{pz1} that, for all $t>0$ and $u\in
L^2_\vt$, $\ck(t,u)$ has a unique extension to a Hilbert-Schmidt
operator from $U^*$ into $L^2_\vt$.
Thus extended, $\ck(t,\punt)$ becomes a bounded linear
operator from $L^2_\vt$ into $L_2(U^*,L^2_\vt)$. Therefore, if $u$
is an $L^2_\vt$-valued adapted process, we can define the stochastic integral as follows:
\beq \int_0^t \left(\Gam(t-s)* \sig(u(s))\right)\, d\cw^*_s := \int_0^t \ck(t-s,\sig(u(s)))\, d\cw^*_s.
\label{int-pz}
\eeq
In the formulation above, $\sig(u(s))$ denotes the function $\sig(u(s))(x):= \sig(u(s,x))$, $x\in \red$, which belongs to $L^2_\vt$.

This definition of the stochastic integral (\ref{int-pz}) is the one
that is used in the mild formulation (\ref{8}).  The main result in
\cite{pz2} on existence and uniqueness of a solution to equation
(\ref{8}) is the following (see \cite[Theorem 0.1]{pz2}).

\begin{theorem}\label{thm-pz}
Assume that $d\in\{1,2,3\}$ and that the coefficients $b$ and $\sig$ are Lipschitz functions. Suppose that there is $\kappa>0$ such that $\Lambda + \kappa\, dx$ is a nonnegative measure (where $dx$ denotes Lebesgue measure), and the spectral measure $\mu$ satisfies
\beq
\int_{\red} \frac{\mu(d\xi)}{1+|\xi|^2}<+\infty.
\label{thm-pz1}
\eeq
Then, for arbitrary $u_0\in H^1_\vt$ and $v_0\in L^2_\vt$, there exists a unique $L^2_\vt$-valued solution to equation (\ref{8}).
\end{theorem}

  As we mentioned at the beginning of this section, this result in \cite{pz2} was extended in \cite{peszat} to higher spatial dimensions.

\subsection{Relation with the random field approach}
\label{relation}

   We now examine the relationship between the random field solution to equation \eqref{2.2} and the $L^2_\vt$-valued solution to equation \eqref{8}.
For this, we assume that the cylindrical Wiener process $W$ considered in the beginning of Section \ref{rva} and the cylindrical $Q$-Wiener process $\cw^*$ (with $Q = \mbox{Id}_{U^*}$) that appears in \eqref{8} are related as follows.

   Let $(e_j)_j$ be a complete orthonormal basis of the Hilbert space $U$ such that $e_j\in \tilde \cs(\red)$, for all $j\geq 1$. Assume that the $e_j^*$ and the $\beta_j(t)$ that appear in \eqref{wstarrd} are given by
\beq\label{2.8a}
   e_j^* = I(e_j)\qquad\mbox{and}\qquad \beta_j(t) = W_t(e_j),
\eeq
where $I$ is the isometry described in Remark \ref{r2}. Recall that $J^*:U^*\rightarrow U_1^*$ denotes a Hilbert-Schmidt embedding between $U^*$ and a possibly larger Hilbert space $U_1^*$; moreover, by Proposition \ref{p1}, $(J^*(e_j^*))_j$ defines a basis in $U_1^*$.

Let us consider $\{u(t,x),\ (t,x)\in [0,T]\times \red\}$, the mild random field solution of (\ref{2.2}) as given in Theorem \ref{existence}, in the case where $L$ is the wave operator in spatial dimension $d=1,2$ or $3$ (so as to have a specific form for $I_0(t,x)$). Then for all $(t,x)\in [0,T]\times \red$,
\begin{align}
u(t,x)=&  \int_{\red} v_0(x-y)\, \Gam(t,dy) + \frac{\partial}{\partial t} \left(\int_{\red} u_0(x-y) \, \Gam(t,dy)\right)\nonumber \\
& + \int_0^t \int_{\mathbb{R}^{d}} \Gamma(t-s,x-y) \sigma(u(s,y))\, W(ds,dy) \nonumber\\
&+ \int_0^t \int_{\mathbb{R}^{d}} b(u(t-s,x-y))\, \Gamma(s,dy) \, ds, \qquad a.s.,
\label{2.9}
\end{align}
The initial conditions $u_0$ and $v_0$ satisfy the hypotheses specified in Lemma \ref{lemma2}. The coefficients $\sig$ and $b$ are Lipschitz functions. Recall that $(t,x) \mapsto u(t,x)$ is mean-square continuous and satisfies
\beq
\sup_{(t,x)\in [0,T]\times\red} E(|u(t,x)|^2)<+\infty.
\label{54}
\eeq

This section is devoted to proving the following result.

\begin{proposition}\label{rdequiv}
Let $\{u(t,x),\, (t,x)\in [0,T]\times \red\}$ be the mild random field solution of \eqref{2.9}. Let $u(t) = u(t,\star )$. Then $\{u(t),\ t \in [0,T]\}$ is the mild $L^2_\vt$-valued solution of \eqref{8}.
\end{proposition}

\noindent\textit{Proof.}
In view of the integral equations (\ref{2.9}) and (\ref{8}), it is clear that the most delicate part in the proof corresponds to the analysis of the stochastic integral terms. Hence, we will start by assuming that both the initial conditions and the drift term $b$ vanish. In this case, $\{u(t,x),\, (t,x)\in [0,T]\times \red\}$ solves the integral equation
$$u(t,x)=  \int_0^t \int_{\mathbb{R}^{d}} \Gamma(t-s,x-y) \sigma(u(s,y))\, W(ds,dy), \quad a.s.$$
for all $(t,x)\in [0,T]\times \red$.
Let us use the following notation for the above stochastic integral:
$$\ci(t,x):=\int_0^t \int_{\red} \Gam(t-s,x-y) \sig(u(s,y))\, W(ds,dy).$$
For any $(t,x)\in [0,T]\times \red$, the above integral is a real-valued random variable and it is well-defined because the integrand satisfies the hypotheses described in Section \ref{distri}, that is, $\Gam(t-\punt,x-\star )$ verifies Hypothesis \ref{hypA} and $\{\sig(u(s,y)),\ (s,y)\in [0,t]\times \red\}$ is a predictable process such that
\beq
\sup_{(s,y)\in [0,T]\times \red} E(|\sig(u(s,y))|^2)\leq C\left( 1+ \sup_{(s,y)\in [0,t]\times \red} E(|u(s,y)|^2)\right)<+\infty.
\label{14.4}
\eeq

Let $u(t) = u(t,\star )$, $t\in [0,T]$. We aim to prove that $\{u(t),\ t \in [0,T]\}$ defines a square-integrable stochastic process with values in the weighted space $L^2_\vt$ which satisfies
$$u(t)=\int_0^t \Gam(t-s)* \sig(u(s))\, d\cw^*_s,\qquad t\in [0,T].$$
Hence, our objective is to prove that $\{\ci(t,\star),\; t\in [0,T]\}$ defines an element in $L^2(\Om\times [0,T];L^2_\vt)$ and
$$\ci(t,\star )=\int_0^t \Gam(t-s)* \sig(u(s))\, d\cw^*_s.$$

 In order to simplify the notation, we will write $Z(s,y):= \sig(u(s,y))$ and let $Z(s)$ denote the function $Z(s)(y)= Z(s,y)$, $y\in\red$.

We will split the proof into several steps.\\

\noindent {\it{Step 1.}} We shall check that  $\{\ci(t,\star ),\; t\in [0,T]\}$ belongs to $L^2(\Om\times [0,T];L^2_\vt)$ and that, for any fixed $(t,x)\in [0,T]\times \red$, the real-valued stochastic integral $\ci(t,x)$ can be written as a stochastic integral with respect to a Hilbert-space-valued Wiener process.

Notice that the norm of $\ci(\cdot,\star )$ in  $L^2(\Om\times [0,T];L^2_\vt)$ coincides with the norm of $u(\cdot)$ in the same space, and the latter is given by
$$E\left( \int_0^T dt \int_{\red} dx\;  \vt(x)\, |u(t,x)|^2 \right).$$
By (\ref{54}) and the fact that $\vt$ is integrable over $\red$, this quantity is finite. In particular, we also deduce that $Z$ belongs to
$L^2(\Om\times [0,T];L^2_\vt)$.

On the other hand, let us recall that $\ci(t,x)$ is a stochastic integral with respect to the cylindrical Wiener process $\{W_s(h),\; s\in[0,T],\, h\in U\}$ (see Section \ref{prel}) with covariance operator $Q=\mbox{Id}_U$ and $s \mapsto \Gamma(t-s,x-\star )Z(s)$ is a predictable process in $L^2(\Omega\times[0,T], U)$ by Proposition \ref{prop1}. Hence, by Proposition \ref{prop24}, the stochastic integral $\ci(t,x)$ may be written as
\beq
\ci(t,x)=\int_0^t \Phi_s^{t,x}\, d\cw_s,
\label{72}
\eeq
where similar to \eqref{QWP},
$$
   \cw_t = \sum_{j=1}^\infty W_t(e_j) J(e_j),
$$
$J:U \to U_1$ is a Hilbert-Schmidt embedding from $U$ into a possibly larger space $U_1$ (note that $U_1$ need not be the dual of $U_1^*$ mentioned after \eqref{2.8a}), and $\{\Phi_s^{t,x}, s\in [0,t]\}$ is the predictable and square integrable process with values in the space $L_2(U,\re)$ of Hilbert-Schmidt operators from $U$ into $\re$, given by
$$\Phi_{s}^{t,x} (h) = \langle \Gam(t-s,x-\star ) Z(s), h \rangle_U,\qquad h \in U.$$
Moreover,
$$E\left(\int_0^t \|\Phi_s^{t,x}\|^2_{L_2(U,\re)} \; ds\right)= E
\left(\int_0^t \|\Gam(t-s,x-\star ) Z(s)\|^2_U \; ds\right).$$

\noindent {\it{Step 2.}} Recall that we aim to prove that
\beq
   \ci(t,\star )=\int_0^t (\Gam(t-s)*Z(s))\, d\cw^*, \quad t\in [0,T],
\label{60}
\eeq
where this equality must be understood in $L^2(\Om\times [0,T]; L^2_\vt)$.

Let $t\in [0,T]$ and $(f_k)_k$ be a complete orthonormal basis in $L^2_\vt$. We will find a suitable expansion of $\ci(t,\star )$ in terms
of $(f_k)_k$. Indeed, by (\ref{72}) and since $\ci(t,\star )$ defines a square integrable $L^2_\vt$-valued random variable, we have the following representation:
\beq
\ci(t,\star )=\sum_{k=1}^\infty \left[\int_{\red} dx \, \vt(x)  \left(
\int_0^t \Phi_s^{t,x} d\cw_s\right)\cdot f_k(x) \right] f_k.
\label{58}
\eeq
Then, by definition of the
stochastic integral with respect to $\cw$ and using representation (\ref{3}) in Proposition
\ref{prop3} (for $H=\re$), for all $x\in \red$,
\begin{align}
\int_0^t \Phi_s^{t,x}\, d\cw_s & = \int_0^t \Phi_s^{t,x} \circ J^{-1}\, d\cw_s \nonumber \\
& = \sum_{j=1}^\infty \int_0^t \Phi^{t,x}_s \circ J^{-1} ( J(e_j))\, dW_s(e_j), \nonumber\\
& = \sum_{j=1}^\infty \int_0^t \Phi^{t,x}_s (e_j)\, d\beta_j(s),
\label{57}
\end{align}
where we have made use of \eqref{2.8a}. Hence, plugging
(\ref{57}) into (\ref{58}), we see that
\beq
\ci(t,\star )=\sum_{k=1}^\infty \left[\int_{\red} dx \, \vt(x)  \left(
\sum_{j=1}^\infty \int_0^t \Phi^{t,x}_s (e_j) \, d\beta_j(s)\right)\cdot
f_k(x)\right] f_k.
\label{59}
\eeq

\noindent {\it{Step 3.}}
We now give an analogous representation for the stochastic
integral on the right-hand side of (\ref{60}). For this,
we will again apply Proposition \ref{prop3} directly to the right-hand side of (\ref{60}); notice that here, $H=L^2_\vt$, and the $J^*$ in \eqref{wstarrd} cancels with the $(J^*)^{-1}$ in the definition of the stochastic integral. Therefore, taking \eqref{2.8a} into account, we see that
\begin{align}
& \int_0^t \Gam(t-s)*Z(s)\, d\cw^* \nonumber \\
& \qquad = \sum_{k=1}^\infty \left(
\sum_{j=1}^\infty \int_0^t \langle \Gam(t-s)*\left(
Z(s)I(e_j)\right),f_k\rangle_{L^2_\vt} \; d\beta_j(s)\right) f_k,
\label{61}
\end{align}
where $(f_k)_k$ and $\beta_j$ are as in Step 2.
Recall that, on the left-hand side of (\ref{61}), $\Gam(t-s)*Z(s)$
is the formal notation for the Hilbert-Schmidt operator defined on
$U^*$ and taking values in $L^2_\vt$ such that, for any $\eta\in
U^{0,*}$, $(\Gam(t-s)*Z(s))(\eta)=\ck(s,Z(s))(\eta)=
\Gam(t-s)*\left( Z(s)\eta \right)$.

  By Lemma \ref{lemma3}, $I(e_j)=e_j*f$ (because $e_j\in \tilde \cs(\red)$), so
equality (\ref{61}) can be written in the form
\begin{align*}
& \int_0^t \Gam(t-s)*Z(s) \, d\cw^*  \\
& \; = \sum_{k=1}^\infty \left( \sum_{j=1}^\infty \int_0^t
\left( \int_{\red} dx\, \vt(x)\; \left[\Gam(t-s)*
\left(Z(s)(e_j*f)\right)\right](x) \cdot f_k(x)\right)
d\beta_j(s)\right) f_k. 
\end{align*}
Applying Fubini's Theorem and comparing the latter expression with (\ref{59}), we observe that, in
order to prove (\ref{60}), it suffices to check that, for almost all $x\in \red$ and any $\ffi\in \tilde
\cs(\red)$,
$$\Phi_s^{t,x}(\ffi)=\left[\Gam(t-s)* \left(Z(s)(\ffi*f)\right)\right](x),\quad s\in [0,t].$$
By definition of the operator $\Phi_s^{t,x}$ and expanding the convolutions on the right-hand side above, this equality is equivalent to
$$\langle \Gam(t-s,x-\star )Z(s),\ffi\rangle_U = \int_{\red} \Gam(t-s,dz) Z(s,x-z) \int_{\red} dy \;
f(x-z-y) \ffi(y).$$ Notice that this is precisely the statement of Lemma \ref{rf} below. Therefore, we can conclude that
(\ref{60}) holds.\\

\noindent {\it{Step 4.}} Let us finally sketch the extension of what we have proved so far to the case of equations (\ref{2.9}) and (\ref{8}). That is, we consider a general Lipschitz continuous drift $b$ and initial conditions $u_0, v_0$ satisfying the hypotheses specified at the beginning of the section. Hence, $\{u(t,x),\; (t,x)\in [0,T]\times \red\}$ satisfies \eqref{2.9}.

   One proves that the process $\{u(t),\; t\in [0,T]\}$ belongs to $L^2(\Om\times [0,T];L^2_\vt)$ as we have done in Step
1. Indeed, an immediate consequence of the proof of Theorem \ref{existence} is that each term in equation (\ref{2.9})
is bounded in $L^2(\Om)$, uniformly with respect to $(t,x)\in [0,T]\times \red$. This clearly implies that each term in
(\ref{2.9}) defines an element in $L^2(\Om\times [0,T];L^2_\vt)$.

   It follows that the stochastic integral $\int_0^t \Gam(t-s)*\sig(u(s))\, d\cw^*$ is well-defined and, by Steps 2 and 3 above, we have
$$\int_0^t \Gam(t-s)*\sig(u(s))\, d\cw^* = \int_0^t\int_{\red}\Gam(t-s,\star -y) \sig(u(s,y)) \, W(ds,dy),$$
where the $\star$ symbol on the right-hand side stands for the variable in $L^2_\vt$.

Concerning the pathwise integral in (\ref{8}), we have
\begin{align*}
\int_0^t \Gam(t-s) * b(u(s))\, ds & = \int_0^t ds\, \int_{\red} \Gam(t-s,dy)\, b(u(s,\star -y))\\
& = \int_0^t ds\, \int_{\red} \Gam(s,dy)\, b(u(t-s,\star -y)).
\end{align*}
It is also clear that the contributions of the initial conditions in equations (\ref{2.9}) and (\ref{8}) coincide as
elements in $L^2([0,T];L^2_\vt)$. We have therefore proved that $\{u(t),\; t\in [0,T]\}$ is the mild solution of (\ref{8}), which concludes the proof of Proposition \ref{rdequiv}.
\qed

\vspace{0.3cm}

We now state and prove the following technical lemma, which was used in the proof of Proposition \ref{rdequiv}.

\begin{lemma}\label{rf}
Fix $t\in [0,T]$. Then, for all $\ffi\in \tilde \cs(\red)$ and $x\in \red$, the stochastic process
$\{\Phi^{t,x}_s(\ffi),\; s\in [0,t]\}$ given by
$$\Phi^{t,x}_s(\ffi)=\langle \Gam(t-s,x-\star )Z(s),\ffi\rangle_U$$
coincides, as an element in $L^2(\Om\times [0,t])$, with $\{\ck^{t,x}_s(\ffi),\; s\in [0,t]\}$, where
$$\ck^{t,x}_s(\ffi)= \int_{\red} \Gam(t-s,dz) Z(s,x-z) \int_{\red} dy \; f(x-z-y) \ffi(y).$$
\end{lemma}

\noindent\textit{Proof.} In order to prove the statement, we will first approximate $\{\Phi^{t,x}_s(\ffi),\; s\in [0,t]\}$ by a
sequence of \emph{smooth} processes.

More precisely, as it has been explained in \cite[Proposition 3.3]{nualartquer}, for any $(s,x)\in [0,t]\times \red$,
we can regularize the element $\Gam(t-s,x-\star ) Z(s)$ of $U$ by means of an approximation of the identity
$(\psi_n)_n\subset C_0^\infty(\red)$, and we can assume that $\psi_n$ is symmetric, for all $n$, and $\vert \tf \psi_n \vert \leq 1$. Then, for any $s\in
[0,t]$, set $J^{t,x}_n(s):= \psi_n * \left( \Gam(t-s,x-\star ) Z(s)\right)$. Again by \cite[Proposition 3.3]{nualartquer},
$J^{t,x}_n(s)$ belongs to $\cs(\red)$ and, as $n\to \infty$, $J^{t,x}_n$ converges to $\Gam(t-\punt,x-\star )
Z$ in $L^2([0,t]\times \Om; U)$. Define
$$\Phi_{n,s}^{t,x} (h) := \langle J^{t,x}_n(s), h \rangle_U,\qquad h \in U.$$
This operator is well-defined because $J^{t,x}_n(s)$ is a smooth function and, in fact, it defines an element in
$L^2\left([0,t]\times \Om; L_2(U,\re)\right)$.

Moreover, $\Phi_n^{t,x} \to \Phi^{t,x}$ in $L^2\left([0,t]\times \Om; L_2(U,\re)\right)$, as $n \to \infty$. Indeed, this is an immediate consequence of the fact that the norm of $\Phi_n^{t,x} -\Phi^{t,x}$ is
given by
\begin{align*}
& E\left( \int_0^t \|\Phi_n^{t,x} -\Phi^{t,x}\|^2_{L_2(U,\re)}\, ds\right) \\
& \qquad = E\left( \int_0^t \sum_{j=1}^\infty
|\langle J^{t,x}_n(s)-\Gam(t-s,x-\star )Z(s), e_j\rangle_U |^2 \, ds \right) \\
& \qquad = E\left( \int_0^t \|J^{t,x}_n(s)-\Gam(t-s,x-\star )Z(s)\|^2_U \, ds\right),
\end{align*}
where $(e_j)_j$ is a complete orthonormal basis in $U$. The last term above tends to zero because, as mentioned before, $J^{t,x}_n \to \Gam(t-\punt,x-\star )Z$ in $L^2([0,t]\times \Om; U)$.

Therefore,  for any $\ffi\in \tilde\cs(\red)$ (in fact, for any $\ffi\in U$), the sequence of real-valued processes $(\Phi^{t,x}_n(\ffi))_n$ converges to $\Phi^{t,x}(\ffi)$ in $L^2(\Om\times [0,t])$. In particular,
$\Phi^{t,x}_n(\ffi)$ converges weakly to $\Phi^{t,x}(\ffi)$, that is, for any $0\leq a<b\leq t$ and $A\in \tf$,
\beq
E\left( 1_A \int_a^b ds \, \Phi^{t,x}_{n,s}(\ffi) \right) \longrightarrow E\left( 1_A \int_a^b ds\,
\Phi^{t,x}_s(\ffi) \right).
\label{63}
\eeq

   We will conclude the proof by checking that the left-hand side of (\ref{63}) also converges to
\beq
E\left( 1_A \int_a^b ds \, \ck^{t,x}_s(\ffi) \right).
\label{63.1}
\eeq
For this, note that, by definition of $\Phi^{t,x}_{n,s}$, the left-hand
side of (\ref{63}) can be written as
$$E\left( 1_A \int_a^b ds \, \langle J^{t,x}_n(s), \ffi \rangle_U \right).$$
Because  $J^{t,x}_n(s)$ and $\vp$ are smooth functions, we can
explicitly compute the inner product in the above expression:
\begin{align*}
& \langle J^{t,x}_n(s),\vp\rangle_U    =
\int_{\red} dy \int_{\red} dz\; J^{t,x}_n(s,y) f(y-z) \vp(z)  \\
&\quad  =  \int_{\red} dy\;  J^{t,x}_n(s,y) \left(f*\vp\right)(y)   \\
& \quad =  \int_{\red} dy \left( \int_{\red} \Gam(t-s,dz)\; \psi_n(y-x+z) Z(s,x-z) \right) \left(f*\vp\right)(y)  \\
& \quad =  \int_{\red} \Gam(t-s,dz) Z(s,x-z) \left( \int_{\red} dy\;
\psi_n(y-x+z) \left(f*\vp\right)(y) \right),
\end{align*}
and so the term on the left-hand side of (\ref{63}) equals

\beq E\left( 1_A \int_a^b ds \int_{\red} \Gam(t-s,dz) Z(s,x-z)
\int_{\red} dy\, \psi_n(x-z-y) (f*\vp)(y) \right).
\label{998} \eeq

   Since $\varphi \in \tilde \cs(\red)$,  the function
$y\mapsto (f*\vp)(y)$ is continuous in $\red$ and
$\lim_{|y|\rightarrow \infty} (f*\vp)(y)=0$.
This implies that, for any $x,z\in \red$,
$$
\lim_{n\rightarrow \infty} \int_{\red} dy\;  \psi_n(x-z-y)
\left(f*\vp\right)(y) = \left(f*\vp\right)(x-z).
$$
Moreover, because $\psi_n$ and $\vp$ belong to $\tilde
\cs(\red)$, we
can apply the definition of the Fourier transform of tempered
distributions:
\begin{align*}
\left| \int_{\red} dy\;  \psi_n(x-z-y) \left(f*\vp\right)(y) \right|
&
=  \left| \int_{\red} \mu(d\xi)\, \tf \psi_n(x-z-\punt)(\xi) \overline{\tf \vp (\xi)}\right| \\
& \leq \int_{\red} \mu(d\xi)\, |\tf \vp (\xi)| < +\infty.
\end{align*}
Thus, in order to apply the Dominated Convergence Theorem in
(\ref{998}), it remains to prove that
$$
   E\left( 1_A \int_a^b ds \int_{\red} \Gam(t-s,dz)\, \vert Z(s,x-z)\vert \right)<+\infty,
$$
and this follows from the hypothesis on $\Gam$ and the process $Z$.
So we have proved that the limit of (\ref{998}), as $n$ goes to
infinity, is
$$
   E \left( 1_A \int_a^b ds  \int_{\red} \Gam(t-s,dz)\, Z(s,x-z)  \int_{\red} dy\, f(x-z-y) \vp(y) \right).
$$
This shows that the left-hand side of (\ref{63}) converges to (\ref{63.1}), which concludes the proof. \qed




\subsection{Relation with the Dalang-Mueller formulation}
\label{relation-dm}

In this section, we examine the relationship between the mild random field solution to equation \eqref{2.9} and the solution introduced by Dalang and Mueller in \cite{dm}, which is based on the $L^2$-valued stochastic integration 
framework that was summarized in Section \ref{ext-dm}. Let $L_\theta^2$ be the space defined in Section \ref{ext-dm}.

   In \cite{dm}, the authors consider solutions to the following stochastic wave equation in $\red$, for any $d\geq 1$:
\beq
\frac{\partial^2 u}{\partial t^2} (t,x)-\Del u(t,x)= \sig(u(t,x)) \dot W(t,x),
\label{100}
\eeq
with initial conditions 
$$u(0,x)=u_0(x), \quad \frac{\partial u}{\partial t}(0,x)=v_0(x), \quad x\in \red,$$
where $u_0, v_0:\red\rightarrow \re$ are appropriate Borel functions. The noise $\dot W(t,x)$ corresponds essentially to the spatially homogeneous Gaussian noise described in Section \ref{prel} (in \cite{dm} a spatial correlation as described in Remark \ref{rem0} was considered). For simplicity, we will restrict to the case 
where the noise is as defined in Section \ref{prel}, with covariance $f$ as in \eqref{1}, and associated spectral measure $\mu$.

   We denote by $H^{-1}(\red)$ the Sobolev space of distributions such that
$$
   \|v\|^2_{H^{-1}(\red)}:= \int_{\red} d\xi \, \frac{1}{1+|\xi|^2}\, |\tf v(\xi)|^2 < +\infty.
$$
According to \cite[Section 5]{dm}, an adapted $L^2_\theta$-valued process $\{u(t,\star ),\; t\in [0,T]\}$ is a {\em mild $L^2_\theta$-valued solution} to (\ref{100}) if
$t\mapsto u(t,\star )$ is mean-square continuous from $[0,T]$ into $L^2_\theta$ and the following $L^2_\theta$-valued stochastic integral equation is satisfied:
\beq
u(t,\star )= \Gam(t) * v_0 + \frac{\partial}{\partial t}\left( \Gam(t) * u_0\right)
+ \int_0^t \int_{\red} \Gam(t-s,\star -y) \sig(u(s,y)) M(ds,dy),
\label{101}
\eeq
where $\Gam$ denotes the fundamental solution of the wave equation in $\red$. The stochastic integral in (\ref{101})
takes values in $L^2_\theta$ and is defined in the final part of Section \ref{ext-dm}. The main result on existence and uniqueness of solutions to equation (\ref{101}) is the following (see \cite[Theorem 13]{dm}).

\begin{theorem}\label{existence-dm}
Assume that the spectral measure $\mu$ satisfies \eqref{18},
$u_0\in L^2(\red)$, $v_0\in H^{-1}(\red)$ and $\sig$ is a Lipschitz function. Then equation (\ref{101}) has a unique mild $L^2_\theta$-valued solution. 
\end{theorem}

In order to be able to compare the solution of the above equation with the mild random field solution to (\ref{2.9}), we consider space dimensions $d\in \{1,2,3\}$ and we set $b=0$. The main result of this section is the following.

\begin{theorem}\label{real-dm} Let $d\in \{1,2,3\}$, and let $\{u(t,x),\, (t,x)\in [0,T]\times \red\}$ be the mild random field solution of \eqref{2.9} (with $b=0$). Let $u(t) = u(t,\star)$. Then $\{u(t),\ t \in [0,T]\}$ is the $L^2_\theta$-valued solution of (\ref{101}).
\end{theorem}

\noindent {\it{Proof.}} 
For simplicity, we assume that the initial conditions vanish (the extension to the general case is straightforward). Recall that $d \in \{1,2,3\}$ and $\{u(t,x),\; (t,x)\in [0,T]\times \red\}$ satisfies the integral equation
\beq
u(t,x)= \ci_{\Gam,Z}(t,x) 
, \quad a.s.
\label{105}
\eeq
for all $(t,x)\in [0,T]\times \red$, where
$$
   \ci_{\Gam,Z}(t,x):=\int_0^t \int_{\red} \Gam(t-s,x-y) Z(s,y) \, W(ds,dy)
$$
and $Z(s,y):=\sig(u(s,y))$. In order to prove that $\{u(t,\star ),\ t \in [0,T]\}$ is the solution of \eqref{101}, we observe that $u(t,\star ) \in L_\theta^2$ a.s., since
$$
   E(\Vert u(t,\star )\Vert_{L_\theta^2}^2) = \int_{\red} E(u(t,x)^2)\, \theta(x)\, dx < \infty
$$
by \eqref{54}. Next, we note that $t \mapsto u(t,\star )$ from $[0,T]$ into $L_\theta^2$ is mean-square continuous, since
$$
   E(\Vert u(t,\star ) - u(s,\star )\Vert_{L_\theta^2}^2) = \int_{\red} E((u(t,x)- u(s,x))^2)\, \theta(x)\, dx,
$$
and we observe that as $s \to t$, by \eqref{54} and since $(t,x) \mapsto u(t,x)$ is $L^2(\Omega)$-continuous, the right-hand side converges to $0$ by the Dominated Convergence Theorem.

   Define 
$$
   v^\theta_{\Gamma,Z}(t,\star ) = \int_0^t \int_{\red} \Gam(t-s,\star -y) \sig(u(s,y))\, M(ds,dy),
$$
where the stochastic integral is defined as in \eqref{103b}. It remains to show that
\beq\label{rd110}
   \ci_{\Gam,Z}(t,\star ) =  v^\theta_{\Gamma,Z}(t,\star )
   \qquad\mbox{in } L^2(\Omega\times\red,dP \times \theta(x) dx).
\eeq

   For this, set $Z_n(s,y)=Z(s,y)1_{[-n,n]^d}(y)$, so that, by definition,
$$
   v^\theta_{\Gamma,Z}(t,\star ) = \lim_{n \to \infty} v^\theta_{\Gamma,Z_n}(t,\star ) \qquad \mbox{in } L^2(\Omega\times\red,dP \times \theta(x) dx),
$$ 
where $v^\theta_{\Gamma,Z_n}(t,\star )$ is defined as in \eqref{102}. By Proposition \ref{integrals}, $v^\theta_{\Gamma,Z_n}(t,\star ) = \ci_{\Gam,Z_n}(t,\star )$ in $L^2(\Omega\times\red,dP \times dx)$, therefore also in $L^2(\Omega\times\red,dP \times \theta(x) dx)$. In order to establish \eqref{rd110}, it suffices to show that
\beq\label{rd4.36}
   E(\Vert \ci_{\Gam,Z}(t,\star ) - \ci_{\Gam,Z_n}(t,\star ) \Vert_{L_\theta^2}^2) \longrightarrow 0 \qquad\mbox{as } n \to \infty.
\eeq

    Set $\Gam_k:=\psi_k*\Gam$, where $\psi_k$ is as in b) of Section \ref{ext-dm}. The expectation in \eqref{rd4.36}  is equal to
\begin{align}
& \int_{\red} dx\, \theta(x)\, E \left(|\ci_{\Gam,Z_n}(t,x)- \ci_{\Gam,Z}(t,x)|^2\right) \nonumber \\
& \qquad = \int_{\red} dx\, \theta(x)\, E \left( \int_0^t ds\, \|\Gam(t-s,x-\star )[Z_n(s,\star )-Z(s,\star )]\|^2_U \right) \nonumber \\
& \qquad \leq C (A_1+A_2),
\label{107}
\end{align}
where $C$ is a positive constant and 
\begin{align*}
   A_1&=  \int_{\red} dx\, \theta(x) E \left( \int_0^t ds\, \|\Gam_k(t-s,x-\star )[Z_n(s,\star )-Z(s,\star )]\|^2_U \right),\\
   A_2 &=  \int_{\red} dx\, \theta(x) \\
&\qquad \qquad \times E \left( \int_0^t ds\, \|[\Gam(t-s,x-\star )-\Gam_k(t-s,x-\star )][Z_n(s,\star )-Z(s,\star )]\|^2_U \right).
\end{align*}
By Proposition \ref{prop1} and the definition of $Z_n$,
$$
A_2\leq C \sup_{(r,y)\in [0,T]\times \red} E(|Z(r,y)|^2) \int_0^T ds \int_{\red} \mu(d\xi) |\tf \Gam(s)(\xi)|^2 |1-\tf \psi_k(\xi)|^2,
$$
and this last expression tends to $0$ when $k\to \infty$. 

Concerning the term $A_1$, since $\Gam_k$, $Z_n$ and $Z$ are functions of the space variable, we smooth $Z_n$ and $Z$ by convolving with $\psi_\ell$, so that we can write the $U$-norm  explicitly, then we use Fatou's Lemma and the fact that $\Gamma_k \geq 0$ to see that
\begin{align}
A_1 & \leq  \int_{\red} dx\, \theta(x) E\left( \int_0^t ds \int_{\red} dy \int_{\red} dz  \, \Gam_k(t-s,x-y) \,\vert Z_n(s,y)-Z(s,y)\vert \right. \nonumber \\
& \qquad  \times f(y-z) \Gam_k(t-s,x-z) \,\vert Z_n(s,z)-Z(s,z)\vert \Big).
\label{109}
\end{align}
Let us prove that, for any fixed $k\geq 1$, the right-hand side of (\ref{109}) converges to $0$ as $n\to \infty$. Indeed, in order to apply the Dominated Convergence Theorem, observe that the integrand in (\ref{109}) converges to $0$ pointwise and its absolute value is bounded, up to some positive constant,  by
$$
   \Gam_k(t-s,x-y)\, |Z(s,y)|\, f(y-z) \Gam_k(t-s,x-z)\, |Z(s,z)|,
$$
which is such that
\begin{align*}
& \int_{\red} dx\, \theta(x) \\
& \; \times E\left( \int_0^t ds \int_{\red} dy \int_{\red} dz  \, \Gam_k(t-s,x-y) |Z(s,y)] f(y-z) \Gam_k(t-s,x-z) |Z(s,z)|\right) \\
& \qquad \leq C \sup_{(r,y)\in [0,T]\times \red} E(|Z(r,y)|^2) \int_0^T ds \int_{\red} \mu(d\xi) |\tf \Gam(s)(\xi)|^2 < +\infty.
\end{align*}
Therefore, taking into account that $A_2 \to 0$ as $k\to \infty$, uniformly with respect to $n$, we deduce that the left-hand side of 
(\ref{107}) tends to $0$ whenever $n\to \infty$. This proves (\ref{rd4.36}), and concludes the proof. \qed



\vskip 1in

\begin{tabular}{ll}
Robert C.~Dalang & Llu\'is Quer-Sardanyons\\
Institut de Math\'ematiques & Departament de Matem\`atiques\\  
Ecole Polytechnique F\'ed\'erale\qquad~ & Universitat Aut\`onoma de Barcelona\\
Station 8 & 08193 Bellaterra (Barcelona) \\
CH-1015 Lausanne & Spain \\
Switzerland & \\
\\
robert.dalang@epfl.ch&  quer@mat.uab.cat
\end{tabular}
\end{document}